\documentclass[11pt]{article}
\usepackage{amsmath, amsthm, amssymb, eucal, stmaryrd, setspace}
\usepackage[all]{xy}
\usepackage{url}

\newcommand{\partentry}[1]{\addtocontents{toc}
                        {\small\bfseries#1\hfill\thepage\par}}

\makeatletter
\def\@part[#1]#2{%
    \ifnum \c@secnumdepth >\m@ne
      \refstepcounter{part}
      \partentry{\protect\makebox[2em][l]{\thepart}#1}
    \else
      \partentry{#1}
    \fi
    {\parindent \z@ \raggedright
     \interlinepenalty \@M
     \normalfont\Large\bfseries\thepart\hspace{1em}#2%
     \markboth{}{}\par}%
    \nobreak
    \vskip 3ex
    \@afterheading}
\def\@spart#1{%
    {\parindent \z@ \raggedright
     \interlinepenalty \@M
     \normalfont\Large\bfseries #1\par}
     \nobreak
     \vskip 3ex
     \@afterheading}
\renewcommand\section{\@startsection{section}{1}{\z@}%
                        {-3.5ex \@plus -1ex \@minus -.2ex}
                        {2ex \@plus.2ex}
                        {\large\bfseries}}
\renewcommand\subsection{
                    \@startsection{subsection}{2}{\z@}
                          {1.75ex \@plus.5ex \@minus.2ex}%
                           {-.4em}
                           {\textbf}}
\def\@seccntformat#1{\@ifundefined{#1@cntformat}%
    {\csname the#1\endcsname\quad}
    {\csname #1@cntformat\endcsname}}
\def\section@cntformat{\thesection.~}
\def\subsection@cntformat{(\thesubsection)\ }
\renewcommand*\l@section{\mdseries\small\@dottedtocline{1}{1.5em}{2em}}
\makeatother
\numberwithin{equation}{subsection}
\newtheorem{thm}[equation]{Theorem}
\newtheorem{cor}[equation]
               {Corollary}
\newtheorem{lem}[equation]
               {Lemma}

\newtheorem{pro}[equation]{Proposition}

\theoremstyle{remark}
\newtheorem{defn}[equation]
               {Definition}
\newtheorem{rem}[equation]
               {Remark}
\newtheorem{exa}[equation]
               {Example}
\newtheorem{notn}[equation]
               {Notation}
\newtheorem{conv}[equation]
               {Convention}

\newcommand{\GG}{\mathbb G}

\newcommand{\NN}{\mathbb N}

\newcommand{\ZZ}{\mathbb Z}
\newcommand{\Ocal}{\mathcal O}

\newcommand{\rest}[1]{{}{\textstyle|}_{#1}}
\newcommand{\size}[1]{{\#(}{#1}{)}}

\newcommand{\Aut}{\operatorname{Aut}}

\newcommand{\Gal}{\operatorname{Gal}}

\newcommand{\im}{\operatorname{im}}
\newcommand{\cha}{\operatorname{char}}

\newcommand{\Pic}{\operatorname{Pic}}

\newcommand{\Spec}{\operatorname{Spec}}

\newcommand{\coker}{\operatorname{coker}}

\newcommand{\textsum}{\textstyle{\sum}}

\providecommand{\abs}[1]{\lvert#1\rvert}
\newcommand{\coa}{\abs}
\newcommand{\sta}{\mathsf}
\newcommand{\stC}{\sta C}
\newcommand{\stL}{\sta L}

\def\pmmu{{\pmb \mu}}

\def\ol{\overline}

\def\wt{\widetilde}

\begin{document}
\title{
\textbf
{Quantitative N\'eron theory for torsion bundles}}
\author{Alessandro Chiodo\thanks
{Financially supported by
the Marie Curie Intra-European
Fellowship within the 6th European Community Framework Programme,
MEIF-CT-2003-501940.}
}
\maketitle

\begin{quote} Let $R$ be a discrete valuation ring
with algebraically closed residue field, and
 consider a smooth curve $C_K$ over the field of fractions $K$.
For any positive integer $r$ prime to the residual characteristic,
we consider the finite $K$-group scheme $\Pic_{C_K}[r]$ of \linebreak $r$-torsion
line bundles on $C_K$. We determine
when there exists a finite $R$-group scheme, which is a model of
$\Pic_{C_K}[r]$ over $R$; in other words,
we establish when the N\'eron model of
$\Pic_{C_K}[r]$ is finite.
To this effect, one needs to
analyse the points of
the N\'eron model over $R$, which,
in general, do not represent
$r$-torsion line bundles
on a semistable reduction of $C_K$.
Instead, we
recast the notion of models on a
stack-theoretic base: there, we find
finite N\'eron models, which represent
$r$-torsion line bundles on a stack-theoretic
semistable reduction of $C_K$.
This allows us to
quantify the lack of finiteness of
the classical N\'eron models and
finally to provide an efficient criterion for it.
\end{quote}
\section{Introduction}
\subsection{The existence of a finite group $R$-scheme, which is an $R$-model of $\Pic_{C_K}[r]$.}
\label{sect:question}
Let $R$ be a discrete valuation ring whose residue field is algebraically closed, we consider
a smooth curve $C_K$ of genus $g\ge 2$ over the field of fractions $K$.
For an integer $r>2$  prime to the residual characteristic,
we consider the \'etale and finite group $K$-scheme $\Pic_{C_K}[r]$ formed by
the $r$-torsion line bundles on $C_K$.

We investigate the existence of a
finite group $R$-scheme, which is an $R$-model of
$\Pic_{C_K}[r]$.
Using a result of Serre \cite{Serre}, Deschamps \cite[Lem.~5.17]{De}
shows that if such an $R$-model exists, then
the minimal regular $R$-model of
$C_K$ is a semistable $R$-curve (the argument is attributed to Raynaud and
uses the fact that $r$ is greater than $2$).

Clearly, this condition is not sufficient, see Example \ref{exa:1irr}.
On the other hand, a result of Lorenzini \cite{Lo} provides a sufficient condition:
a finite $R$-model of $\Pic_{C_K}[r]$
equipped with a group structure exists
if the dual graph of the special fibre of the
minimal regular model of $C_K$ over $R$
is obtained from another graph  by dividing all the
edges in $r$ edges.
This condition is not necessary, as we illustrate in \S\ref{sect:examples},
Examples \ref{exa:wheels} and \ref{exa:rdiv}.

This paper places this issue in a different perspective:
describing
the interplay between the finiteness of
the models of $\Pic_{C_K}[r]$ on a stack-theoretic base
and the existence of stack-theoretic
semistable reductions of $C_K$.
As an application we finally find a necessary and sufficient condition on $C_K$
for the existence of a finite group
scheme extending $\Pic_{C_K}[r]$ over $R$.
Indeed, Corollary \ref{cor:group} shows that
the following statements are equivalent.
\begin{enumerate}
\item There exists a finite group $R$-scheme, which is an $R$-model of $\Pic_{C_K}[r]$.
\item The minimal regular model of $C_K$ over $R$ is semistable
and the dual graph of its special fibre satisfies the following property:
the (signed) number of edges common to any two  circuits is always a multiple of $r$.\footnote{We always count the signed number of common edges
as detailed in \eqref{conv}.}
\end{enumerate}
\subsection{The problem in terms of N\'eron models.}
The question posed in \eqref{sect:question}
is about the N\'eron model of $\Pic_{C_K}[r]$.
Since $r$ is prime to
the residual characteristic, by \cite[7.1/Thm.~1]{BLR}, a finite group $R$-scheme, which is a model of
$\Pic_{C_K}[r]$
is necessarily the N\'eron model of $\Pic_{C_K}[r]$, the
universal $R$-model in the sense
of the N\'eron mapping property.
This model is a canonical $R$-model of $\Pic_{C_K}[r]$
equipped with a group structure naturally induced by
the universal property.
In terms of the N\'eron model an equivalent reformulation of
the question above is: when is $N(\Pic_{C_K}[r])$ finite over $R$?

In order to explore N\'eron models of $\Pic_{C_K}[r]$,
it would be convenient to realise them
in terms of the functor of $r$-torsion line bundles
on a semistable reduction of $C_K$ over $R$.
However, this is not possible in general.
First, it may well happen that such a semistable reduction does not exist
over $R$.
Second, if the semistable reduction
$C_K$ exists but is not of compact type\footnote{It contains nodes whose desingularisation
is connected (nonseparating nodes).}
the functor of $r$-torsion
line bundles is not finite, whereas the N\'eron model may well be finite, see
Example \ref{exa:2irr}.

\subsection{Quantitative theory of N\'eron models via twisted semistable reductions.}
We start by slightly recasting the notion of
N\'eron model. We place it over a new stack-theoretic base
$S[d]$: a proper modification of $S=\Spec R$,
containing the generic point $\Spec K$, and
having a stabiliser $\pmmu_d$ over the special point
(we assume $d$ prime to the residual characteristic).
On $S[d]$, there exists
a N\'eron $d$-model
$N_d(\Pic_{C_K}[r])$: a universal $S[d]$-model
in the sense of the N\'eron mapping property,
see Definition \ref{defn:dner}.
In this perspective, the problem is not only investigating
geometric properties (such as the finiteness of N\'eron models, or the
existence of semistable reductions) on a fixed base, but
also measuring how large $\pmmu_d$ should be
for such properties to hold on a stack-theoretic base.
The question posed in \eqref{sect:question} becomes: for which values of $d$
the N\'eron $d$-model of $\Pic_{C_K}[r]$ is finite?

The advantage is that semistable reductions can be more efficiently
used in this setting. Indeed, when $d$ is sufficiently large
we can realise N\'eron models of $\Pic_{C_K}[r]$
via stack-theoretic semistable reductions; then, by descent,
we can more easily understand the picture for lower values of $d$.
We use stack-theoretic curves whose nontrivial stabilisers occur only
at the nodes (twisted curves in the sense of Abramovich and Vistoli).
 As shown in \cite[Thm.~3.2.2]{Ch_mod},
 these curves
 carry as many $r$-torsion
 line bundles as smooth curves as soon as the size of the
 stabilisers on
 all nonseparating
 nodes divides $r$ (we recall the result in
 Theorem \ref{thm:cond}).
Therefore, the problem of constructing finite N\'eron $d$-model is
solved as soon as a twisted semistable reduction of this sort exists.
The study of geometric invariants of $C_K$ allows us to
determine when this is the case.

\subsection{The answer in terms of geometric invariants of $C_K$.}
The following natural invariants of $C_K$ play a crucial role.
We assume that $\Pic_{C_K}[r]$ is tamely ramified over $K$, otherwise, for any $d$,
the N\'eron $d$-model is not finite. Then, we consider

\smallskip
\noindent \begin{tabular}
{lllllll}
--& $m_1$,& the least integer such that the minimal regular model of $C_K$ over $S[d]$ is semistable,
\quad \quad\quad \quad\quad \quad\\
--& $\Gamma$, &  the dual graph of the special fibre of the minimal regular model over $S[m_1]$.
\end{tabular}

\medskip

Using these invariants,
In Theorem \ref{thm:group} and Proposition \ref{pro:m2m1}, we provide explicit formulae
for $m_2$, the least integer for which
the N\'eron $d$-model is finite.
On $S[m_1r]$, there exists a twisted
semistable reduction $\stC$ of $C_K$
with stabilisers of order $r$ on all nonseparating nodes.
As illustrated above, this yields a finite N\'eron $rm_1$-model of
$\Pic_{C_K}[r]$, which represents the
$r$-torsion line bundles over a
twisted curve $\stC$. For $d\mid r$,
we can descend from $S[rm_1]$ to $S[dm_1]$ if
the special fibre of $\stC$ and all its $r$-torsion line bundles
are fixed by the natural action of $\xi_{r}^{d}$.
In this way, the study of the action
on the twisted curve finally yields \eqref{eq:m1m2}
$$m_2= m_1r/\gcd\{r,c\},$$
where $c$ is the greatest common divisor of the
number  of
edges shared by two  circuits in the dual graph $\Gamma$.
In particular, we find the answer to the question posed in \eqref{sect:question}, see Corollary
\ref{cor:group}.

A similar analysis allows us to determine $m_3$, the
least integer for which
the N\'eron $d$-model is finite and represents the $r$-torsion
on a twisted semistable reduction of $C_K$. We have (equation \eqref{eq:m1m3})
$$m_3= m_1r/\gcd\{r,t\},$$
where $t$ is the greatest common divisor of the
thicknesses (see Notation \ref{notn:thickness})
of the nonseparating nodes of the special fibre of
the stable reduction of $C_K$ over $S[m_1]$. In \eqref{sect:examples},
we compare the values of $m_3$ and $m_2$ in several examples.

\subsection{The group of
connected components of the special fibre of $N(\Pic^0_{C_K})$.}
The fact that the N\'eron model of the $\Pic_{C_K}[r]$ need
not represent the functor of $r$-torsion line bundles is well known
and holds in general for the entire functor $\Pic^0_{C_K}$.
In the literature, this issue motivated the introduction of
a finite abelian group $\Phi$:
the group of
connected components of the special fibre of the N\'eron model
of $\Pic^0_{C_K}$. In \cite[\S8.1.2, p. 64]{Re_Pic},
Raynaud
showed how,
when the minimal regular model $C_R$ of $C_K$ is semistable,
the group  $\Phi$ can be defined in terms of
the dual graph $\Gamma$
of the special fibre of $C_R$.
The group $\Phi$
is
studied from very different viewpoints (see \S\ref{sect:literature}),
it is difficult to determine in general, and
it is the subject of
several open questions, see \cite{Lor:preprint} for  recent results.

Using $\Phi$, one can state an evident reformulation of our main question.
The special fibre of
the N\'eron model of $\Pic_{C_K}[r]$
can be regarded as an extension of $\Phi[r]$ by the
group of $r$-torsion line bundles on the special fibre of $C_R$.
In this way, the finiteness of $N(\Pic_{C_K}[r])$ is equivalent to
\begin{equation}\label{eq:H}
\Phi[r]\cong (\ZZ/r\ZZ)^{\oplus b_1},
\end{equation}
where $b_1$ is the first Betti number of the dual graph $\Gamma$.
This is precisely the statement for which
Lorenzini provides the sufficient
condition \cite[Prop.~2]{Lo} mentioned above.
In this respect, our criterion stated in \eqref{sect:question}
can be also regarded as a statement on $\Phi$:
it unravels the geometric condition
on $C_K$ encoded in \eqref{eq:H}.

We find it
interesting to give the proof of
our claim
both from the geometric
point of view of twisted curves (Corollary \ref{cor:group}) and
from the abstract point of view of the group
$\Phi$ (Propositions \ref{pro:viane} and \ref{pro:tor}). Indeed, the approach via line bundles on
twisted curves of the present paper
is a new geometric interpretation of
the $r$-torsion of the group $\Phi$. (Other
interactions between the geometry of curves and
the group $\Phi$ have been recently shown by Caporaso in \cite{Cap}.)
In Section \ref{sect:without} we do not rely
on the theory of twisted curves---we wrote the
entire section  so
that the reader can start directly there if he wishes.

\subsection{The torsor of $r$th roots of a line bundle.}
We also consider the case of $r$th roots of a
given line bundle $F_K$ on $C_K$,
whose degree is a multiple of $r$. The
$r$th tensor power $L\mapsto L^{\otimes r}$
in the Picard group $\Pic_{C_K}$ is \'etale and
$F_K$ is in the image; therefore
the $r$th roots of $F_K$ form an \'etale, finite $K$-torsor $T_K$ over
$K$ under the kernel $G_K=\Pic_{C_K}[r]$.

In Section \ref{sect:torsor}, we generalise
the results obtained for $\Pic_{C_K}[r]$: we find
N\'eron $d$-models of
$T_K$ which represent the $r$th roots of a line bundle extending
$F_K$ over a twisted semistable reduction of $C_K$ over $S[d]$.
The criterion for the finiteness of $N(T_K)$ over $R$
is given in Corollary \ref{cor:torsor}.

\subsection{Structure of the paper.}
In Section \ref{sect:term} we fix the terminology.
In Section \ref{sect:twisted}, we recall some known results on
twisted curves.
In Section \ref{sect:dmodels}, we set up the notion of N\'eron $d$-models.
In Section \ref{sect:group}, we prove the results
mentioned above, Theorem \ref{thm:group} and Corollary \ref{cor:group}.
In Section \ref{sect:torsor}, we generalise the results to the torsor $T_K$.
Section \ref{sect:without} is self contained: we focus on
the criterion stated above, and we prove it by means of the group $\Phi$ and with no use of
the theory of twisted curves.
In Section \ref{sect:circuits} we show a technical lemma.
\subsection{Acknowledgements.}
I would like to thank
Andr\'e Hirschowitz for encouraging me to investigate this problem.
I'm indebted to Michel Raynaud for hours of explanations and
suggestions which greatly influenced this work.
I would like to thank Dino Lorenzini for attentively reading my
early attempts to state Corollary \ref{cor:group} and for
pointing out several results in the existing literature to me.

\section{Terminology}\label{sect:term}
\subsection{Context.}
We work with schemes locally of finite type over an algebraically closed field $k$.
By $R$ we denote a complete discrete valuation ring with residue field $k$ and
fraction field $K$.
We write $r$ for a positive integer prime to $\cha(k)$.

\subsection{Graphs.}
In this paper we say graphs for 2-graphs possibly having loops and
multiple edges such as in

\begin{equation*}
\xymatrix@=1.2pc{
*{ \,\bullet} \ar@{-}@/_/[r]\ar@{-}@/^/[r]&
*{ \,\bullet} \ar@{-}[r]
&*{ \,\bullet} \ar@{-}@(ur,dr) & {}\\ &
}
\end{equation*}
In fact, all graphs considered in this paper arise in the following way.

\begin{notn}[dual graph]\label{notn:dualgraph}
The \emph{dual graph} $\Gamma_C$
of a semistable curve $C$ over $k$ is the graph whose
set of vertices
$V$ is the set of irreducible components
of $C$ and whose set of edges
$E$ is the set of nodes of $C$. An edge $e$ associated to a node in $C$
joins the vertices $v'$ and $v''$
if the following condition is satisfied:
each irreducible component corresponding to $v'$ and
 $v''$ contains a
branch of the node.

We say that an edge $e$ is \emph{nonseparating} if the graph remains connected
after taking off $e$ from $E$.
In this way the nodes of $C$ are either separating or nonseparating according to
the corresponding edge in the dual graph.
\end{notn}
\begin{notn}[paths and circuits]\label{notn:paths}
Each edge in a graph is \emph{oriented} if it is equipped with
an ordering for its vertices: we refer to the first vertex and the last
vertex as the tail and the tip.
A \emph{path} in a graph is a sequence of oriented
edges $e_0, \dots, e_{n-1}$
such that the tip of $e_i$ is also the tail of $e_{i+1}$.
In this way, a path determines a sequence of vertices
$v_0, \dots, v_n$ such that from the vertex $v_i$ there is an
edge $e_i$ to the vertex $v_{i+1}$
for $0\le i<n$.
A \emph{simple path} is a path
whose vertices $v_0,\dots, v_n$ are distinct.
A \emph{circuit}, is a path whose
vertices $v_0$ and $v_n$ coincide and with no
repeated vertices aside from the first and the last.
Note that the edges of paths and circuits are always oriented: these
paths are sometimes called directed; here we omit the adjective
directed, because no ambiguity may arise.
\end{notn}
\begin{conv}[signed number of edges]\label{conv}
When counting the edges common to two circuits in a graph,
we always count the signed number of edges. More precisely,
the edges lying on both circuits are counted with
a positive or negative sign according to whether their orientations
in the two circuits agree or disagree.
\end{conv}
\begin{notn}[chain and cochain complexes]
\label{notn:dualchain}
Let $\Gamma$ be a graph with an arbitrary fixed orientation.
We have a chain complex with differential
\begin{align*}
\partial \colon C_1(\Gamma,\ZZ)&\xrightarrow{\ } C_0(\Gamma,\ZZ),
\end{align*}
where the edge $e$ is sent to $[v_+]-[v_-]$ if $v_+$ is the tip and $v_-$ is the tail of $e$.
Since $C_0(\Gamma,\ZZ)$ and $C_1(\Gamma,\ZZ)$ are canonically
isomorphic to the group of cochains $C^0(\Gamma,\ZZ)$ and $C^1(\Gamma,\ZZ)$ we can regard the
differential of the cochain complex
$C^\bullet(\Gamma,\ZZ)$
as
\begin{align*}
\delta \colon C_0(\Gamma,\ZZ)&\xrightarrow{\ } C_1(\Gamma,\ZZ),
\end{align*}
where the vertex $v$ is sent to the sum of the edges ending at $v$
minus the sum of the edges starting at $v$.

For any positive integer $q$, after tensoring with $\ZZ/q\ZZ$, we get
differentials of $\ZZ/q\ZZ$-valued chain and cochain complexes
\begin{equation}
\label{eq:partialq}
\partial _q\colon C_1(\Gamma,\ZZ/q\ZZ)\xrightarrow{\ } C_0(\Gamma,\ZZ/q\ZZ),\quad \quad
\delta _q\colon C_0(\Gamma,\ZZ/q\ZZ)\xrightarrow{\ } C_1(\Gamma,\ZZ/q\ZZ).
\end{equation}
\end{notn}
\subsection{Stacks.}\label{sect:stacks}
We always
work with algebraic stacks locally of finite type over a locally
noetherian scheme (we refer to \cite{LM} for the main definitions).
When working with algebraic stacks with finite diagonal,
we use Keel and Mori's Theorem
\cite{KM}: there exists an algebraic space $\coa{\sta X}$ (the \emph{coarse space})
associated to
the stack $\sta X$
and a morphism $\pi_{\sta X}\colon \sta X\to \coa{\sta X}$
(or simply $\pi$)
which is universal with respect
to morphisms from $\sta X$ to algebraic spaces.

A \emph{geometric point} $\sta p\in \sta X$
is an object  $\Spec k\to \sta X$.
We denote by $\Aut(\sta p)$ the \emph{stabiliser} of $\sta p$: the
automorphism group
of $\sta p$ regarded as an object of the fibred category $\sta X_k$.
In order to identify certain stacks and morphisms between stacks locally at a geometric point
$\sta p$, we
adopt the standard convention (see \cite[\S1.5]{ACV})
of exhibiting the strict henselisation of the stack and of the morphisms
involved; we call this process ``local picture at $\sta p$'' (this avoids repeated mention
of strict henselisation of the morphisms and of the stacks).

\subsection{Semistable curves.}\label{sect:sstable}
All curves appearing in this paper are semistable:
for sake of clarity we recall the definition.
\begin{defn}\label{defn:sstable}
A \emph{semistable} curve
of genus $g\ge 2$ on
a scheme $X$ is
a proper and flat morphism
$C\to X$
whose fibres $C_x$ over geometric points $x\in X$ are reduced, connected,
$1$-dimensional, and satisfy the following conditions:
\begin{enumerate}
\item $C_x$ has only ordinary double points (the nodes),
\item if $E$ is a nonsingular rational component of $C_x$, then
$E$ meets the other components of $C_x$ in at least two points.
\item we have $\dim_{k(x)}H^1(C_x,\Ocal_{C_x})=g$.
\end{enumerate}
The curve is \emph{stable} if, in (2), we require that
$E$ meets the other components of $C_x$ in at least three points.
\end{defn}
\begin{rem}[line bundles on semistable curves]\label{rem:seqgl}
If $C$ is  a  semistable curve over $k$, an essentially complete description of
$\Pic_C=H^1(C,\GG_m)$
and of its $r$-torsion subgroup $\Pic_C[r]=H^1(C,\pmmu_r)$ is given as follows.
Consider the short exact sequences
$1\to \GG_m \to \nu_*\GG_m\to \nu_*\GG_m/\GG_m\to 1$
and
$1\to \pmmu_r\to \nu_*\pmmu_r\to \nu_*\pmmu_r/\pmmu_r\to 1$.
Let us choose an orientation for the dual graph of $C$. We denote by $\nu\colon C^\nu \to C$ the
normalisation of $C$.
Then, the long exact sequences of cohomology induced by the above sequences are
\begin{eqnarray*}
  &1\to \GG_m\to (\GG_m)^V \to (\GG_m)^E\to \Pic_{C}\xrightarrow{\ \nu^* \ } \Pic_{C^\nu}\to 1,
\end{eqnarray*}
and
\begin{eqnarray*}
  &1\to \pmmu_r\to (\pmmu_r)^V \to (\pmmu_r)^E\to \Pic_{C}[r]\xrightarrow{\ \nu^* \ } \Pic_{C^\nu}[r]\to 1,
\end{eqnarray*}
where $V$ and $E$ are the vertices and the edges of the dual graph and the morphisms
$(\GG_m)^V \to (\GG_m)^E$ and $(\pmmu_r)^V \to (\pmmu_r)^E$ can be regarded as
coboundary homomorphisms $\delta$ with the assigned orientation as in Notation \eqref{notn:dualchain}.
In this way, the number of $r$-torsion line bundles on a semistable curve
coincides with $r^{2g-b_1}$, where $b_1$ is the first Betti number of the dual graph.
\end{rem}
\begin{notn}[thickness]\label{notn:thickness}
Let $C_R$ be a semistable curve over
the complete discrete valuation ring $R$. When the fibre over the fraction field $K$ is smooth,
the local picture at a node $e\in C$ is
given by $\Spec R[z,w]/(zw-\pi^{\eta (e)})$, for $\pi$ a uniformiser
of $R$ and $\eta(e)$ a positive integer, which we call the \emph{thickness} of the node $e$.
\end{notn}
\begin{notn}[stable model and semistable minimal regular model]\label{notn:ssmodels}
Let us assume that a smooth curve $C_K$ of genus $g\ge 2$ over $K$ is given and
that there exists a \emph{semistable reduction}
$C_R$ over $R$, i.e.
$C_R$ is a semistable curve
over $R$ and its generic fibre is isomorphic to $C_K$.
In this situation, there may be
several semistable reductions, but among them there are
two special choices.
Indeed, there exists a unique stable curve
$C_R^{\rm st}\to \Spec R$ whose generic fibre is isomorphic to $C_K$;
we refer to it
as the \emph{stable model of $C_K\to \Spec K$}.
Furthermore,  there exists a {unique} semistable curve
$C_R^{\rm reg}\to \Spec R$ for which
$C_R^{\rm reg}$ is regular and
we refer to it as the
\emph{semistable minimal regular model of $C_K\to \Spec K$}
(indeed $C^{\rm reg}_R$ is minimal with respect to
regular models of $C_K$).

There is a natural $R$-morphism $C_R^{\rm reg}\to C_R^{\rm
st}$ obtained by contraction of all rational lines of
selfintersection $-2$ in $C_R^{\rm reg}$. Each node of thickness
$\eta(e)$ in $C_R^{\rm st}$ is the contraction of a chain of $\eta(e)-1$
rational curves in $C_R^{\rm reg}$.
The statements above are
well known from the
theory of semistable
reduction \cite{SGA7} \cite{DM}.
\end{notn}
\subsection{Line bundles on semistable reductions: multidegrees and dual graph.}\label{sect:multideg}
Consider the semistable minimal regular model $C_{R}^{\rm reg}$ of $C_K$.
We study line bundles $F$ on $C_R^{\rm reg}$ whose restriction on
$C_K$ is trivial: $F\rest{C_K}\cong \Ocal_{C_K}$.

They can be written as $\Ocal (\sum _{v\in V} a_v X_v)$, where $V$ is the
set of irreducible components $X_v$ of the special fibre and
$\vec{a}=(a_v)_V$ is a multiindex with entries in $\ZZ$. Note that
the multiindex $\vec{d}(F)=(d_v)_V$ given by the degree of
$F=\Ocal(\sum _{v\in V} a_v X_v)$ on all irreducible components
satisfies $$\vec{d}(F)=M\vec{a},$$
where $M$ is the intersection matrix $(X_{v_1}\cdot X_{v_2})$.

This shows that the set of multiindices determined by the degrees
on all irreducible components of all
line bundles $F$ extending $\Ocal_{C_K}$ on $C_R^{\rm reg}$
equals  $\im(M)\subset C_0(\Gamma, \ZZ)$, where $\Gamma$ is the dual graph of the special fibre of
$C_R^{\rm reg}$ and $M$ is regarded as an endomorphism of
$C_0(\Gamma,\ZZ)$.
\begin{rem}\label{rem:Mdelta}
It is easy to see that we have $-M=\partial \circ \delta.$\label{}
\end{rem}
\begin{rem}\label{rem:FonCreg} For any line bundle $F_K$,
whose relative degree is a multiple of $r$, it is easy to
see that the following conditions are equivalent.
\begin{enumerate}
\item There exists a line bundle  $F$ on $C_R^{\rm reg}$ satisfying $F\rest{C_K}\cong F_K$,
whose degree is a multiple of $r$ on each irreducible component of the special fibre.
\item There exists a line bundle  $F$ on $C_R^{\rm reg}$ satisfying $F\rest{C_K}\cong F_K$
and such that the multiindex $\vec{d}(F) \mod r$ belongs to
 $\im (\partial_r \circ \delta_r)$.
\end{enumerate}
\end{rem}
\subsection{Twisted curves.}\label{sect:subtwisted}
Consider a
proper and flat
morphism from a stack $\stC$ to a scheme $X$, for which
the fibres are purely one-dimensional
with at most nodal singularities,
the order of all stabilisers is prime to $\cha(k)$,
the coarse space is a semistable curve
$\coa{\sta C}\to X$ of genus $g$, and
the smooth locus
 is a scheme.

These stack-theoretic curves $\stC\to X$ are called
\emph{twisted curves} if the following condition
introduced by Abramovich and Vistoli \cite{AV} is satisfied.
The local picture
at a node
is given by $[U/\pmmu_l]\to T$, where
\begin{enumerate}
\item[$\bullet$] $T=\Spec A$ for $A$ a ring,
\item[$\bullet$] $U=\Spec A[z,w]/(zw-t)$ for some $t\in A$,
\item[$\bullet$] the action of $\pmmu_l$ is given by
$(z,w)\mapsto (\xi_lz,\xi_l^{-1}w)$ for $l$ a positive integer and $\xi_l$ a primitive $l$th root of unity.
\end{enumerate}

\begin{rem}[unbalanced twisted curves]
In the existing literature
twisted curves satisfying the above
local condition
are often called \emph{balanced}, \cite{AV}.
We drop the adjective balanced, because
we never consider unbalanced twisted curves.
\end{rem}

All the notions introduced above for semistable curves generalise
to twisted curves.
\begin{itemize}
\item[--] Clearly,
the notion of dual graph extends word for word
from semistable curves over $k$ to twisted curves over $k$.

\item[--] Furthermore, the exactness of the sequences
$1\to \GG_m \to \nu_*\GG_m\to \nu_*\GG_m/\GG_m\to 1$
and
$1\to \pmmu_r\to \nu_*\pmmu_r\to \nu_*\pmmu_r/\pmmu_r\to 1$
holds for twisted curves as well as for semistable curves. However note that, in general,
the long exact sequences written in Remark \ref{rem:seqgl} are exact only up to the
last homomorphism in the sequence, the pullback
$\nu^*$, which in general is not surjective. This happens because
the higher cohomology groups do not vanish, see \cite[\S3]{Ch_mod} for more details.

\item[--] For a twisted curve $\stC_R$ over $R$
with smooth generic fibre, the notion of thickness of a node $\sta e$ also extends.
The thickness $\eta(\sta e)$ is a positive integer  such that the local picture
of $\stC_R$ at $\sta e$ is given by $[U/\pmmu_l]$ where $U$ is
the scheme $\Spec R[z,w]/(zw-\pi^{\eta (\sta e)})$, for $\pi$ a uniformiser
of $R$.

\item[--] Given a smooth curve $C_K$ over
the fraction field $K$ of a discrete valuation ring $R$,
we say that \emph{$\stC_R$ is a twisted semistable reduction of $C_K$ on $R$}
if $\stC_R\to \Spec R$ is a twisted curve
over $R$ whose  generic fibre is isomorphic to $C_K$.

\end{itemize}
\begin{rem}[adding cyclic stabilisers on thick nodes]\label{rem:modifsst}
Let $C_K$ be a smooth curve and let $C_R$ be a
semistable reduction. By definition, $C_R$
is a twisted semistable reduction, because a semistable reduction
is just a representable twisted semistable reduction.
Let $e$ be a node of thickness $\eta$ in the special fibre of $C_R$.
Assume that $\eta$ factors as
$\eta=dh$ where $d$ is a prime to $\cha(k)$.
Then, we can construct a twisted semistable reduction
$$C(d)\to \Spec R$$ with a
stabiliser of order $d$ overlying $e$.

Indeed, take an \'etale neighbourhood of $e\in C_R$ of the form $\ol U=\Spec
R[w,z]/(zw=\pi^{\eta})$. Consider the quotient
stack $[U/\pmmu_d]=[\{zw=\pi^{\eta/d}\}/\pmmu_{d}]$ with $\pmmu_{d}$
acting as $(z,w)\mapsto(\xi_{d}z,\xi_{d}^{-1}w)$.
The action is free outside the origin $\sta p=(z=w=0)$, and
$[U/\pmmu_d]\to \ol U$ is invertible away from $\sta p$.
We define a twisted reduction $C_R(d)$ by glueing $[U/\pmmu_d]$ to
$C\setminus \{e\}$ along the isomorphism
$[U/\pmmu_d]\setminus \{p\}\to \ol U\setminus\{e\}$.
Note that $C(d)$ has thickness $\eta/d$ at $\sta p\in [U/\pmmu_d]$.
\end{rem}

\begin{rem}\label{rem:char0}
We point out that this shows that---in characteristic $0$---any
curve $C_K$ admitting a
semistable reduction also admits a twisted semistable reduction
over $R$ which, as a stack,  is  \emph{regular}.
This fact explains why extending line bundles is easier
if we allow twisted
semistable reductions.
\end{rem}

\section{Taking $r$th roots on twisted curves}\label{sect:twisted}
\subsection{Line bundles on twisted curves.}\label{sect:monod}
Consider the following situation: $\sta L$ is a line bundle
on a twisted curve over $k$. Locally at a node
whose stabiliser has order $l$, the line bundle
$\sta L$ can be regarded as the trivial line bundle
over $\{xy=0\}$ with $\pmmu_l$-linearisation
\begin{equation*}\label{eq:stablb}((x,y),\lambda)\mapsto ((\xi_lx,\xi_l^{-1}y),\xi_l^{q}\lambda).\end{equation*}
We point out that the multiplicity index $q$ is
uniquely determined in $\{0,\dots, l-1\}$ whenever
one of the branches of the node  is specified.
Indeed, $q\in \{0,\dots, l-1\}$
is determined
by the line bundle $\sta L$, the node, and the choice of the branch as follows.
The tangent line along the chosen branch can be regarded as the $\pmmu_l$-representation
  $T\colon z\mapsto \xi_l z$ for a suitable primitive $l$th root of unity.
Similarly,
the restriction of $L$ to the chosen branch is a $\pmmu_l$-representation and
is a power of order $q\in \{0,\dots, l-1\}$ of $T$.
\subsection{The functor of $r$th roots.}
For any twisted curve $\stC$ on $S=\Spec R$
and for any line bundle $\sta F$ on $\stC$ whose relative degree is a multiple of $r$, the $r$th roots of
$\sta F$ on $\stC$ are represented by an \'etale $S$-scheme $$\sta F^{1/r}\to S.$$
More precisely, following \cite[\S3]{Ch_mod},
we can consider
the category $\sta {Roots}_r$ formed by the objects $(T, M_T, j_T)$,
where $T$ is an $S$-scheme, $M_T$ is a line bundle on
the semistable curve $C\times _S T\to T$,
and $j_T$ is an isomorphism
of line bundles $M_T^{\otimes r}\xrightarrow{\sim} \sta F\times_S T$.
The category  $\sta {Roots}_r$
is a stack of Deligne--Mumford type, \'etale on $S$,
\cite[Prop.~3.1.3]{Ch_mod}.

In order to obtain the scheme $\sta F^{1/r}$, let us point out that
every object $(T, M_T, j_T)$ has an automorphism
given by multiplication by an $r$th root of unity on $T$
along the fibre of $M_T$. Then, $\sta F^{1/r}\to S$
can be regarded as the ``rigidification along
$\pmmu_r$ of $\sta {Roots}_r$'' in the sense of
Abramovich, Corti, and Vistoli, \cite[\S5]{ACV} (denoted by ``$\!\!\fatslash$ '', see \cite{Ro} for
a careful treatment of this issue)
$$\sta F^{1/r}=\sta {Roots}_r \!\!\fatslash {\pmmu_r}.$$
This process of
rigidification provides a natural framework to a standard procedure that
occurs systematically in the construction of  Picard functors.
Indeed, for any morphism of schemes $f\colon Y\to S$,
the natural functor $T\mapsto \Pic(Y_T)$ from $X$-schemes
to sets is in general a presheaf and is not represented by a scheme.
The actual ``relative Picard functor'' $\Pic_{Y/S}$ is defined
by the passage to
the associated sheaf (see \cite[Ch.~8]{BLR}).
In this way, when $\stC$ is a semistable curve on $S$,
the construction of $\sta F^{1/r}$
yields the same scheme as the subscheme of $r$th roots of $F$ in the
relative Picard functor $\Pic_{C/S}$.
In particular, the $K$-group
$\Pic_{C_K}[r]$ equals $\Ocal ^{1/r}$ and the $K$-torsor  of
$r$th roots of $F_K$ in $\Pic_{C_K}$ equals $F_K^{1/r}$.

The scheme $\sta F^{1/r}$ is \'etale on $S$.
If we assume that $C_K$ is a smooth curve, then
the generic fibre of the scheme $\sta F^{1/r}$ contains $r^{2g}$ points.
In this way, ${\Ocal^{1/r}}$ is a finite group scheme
and ${\sta F^{1/r}}$
is a finite torsor under ${\Ocal^{1/r}}$ if and only if
the restriction $\sta F\rest {\stC_k}$ to the special fibre
$\stC_k$ has exactly $r^{2g}$ $r$th roots.
This motivates the following result.
\begin{thm}[{\cite[Thm.~3.2.2]{Ch_mod}}]\label{thm:cond}
Let $\pi\colon \stC\to \coa{\stC}$ be
a twisted curve of genus $g$ over $k$.
Let $F$ be a line bundle on $\coa{\stC}$, whose total degree is a multiple of $r$.
There are exactly $r^{2g}$ roots of $\pi^*F$ on $\stC$ if and only if
\begin{align*}
&\size{\Aut(\sta e)} \in r\ZZ && \text{for each nonseparating node $\sta e$, and}\\
&\size{\Aut(\sta e)} d(\sta e) \in r\ZZ&& \text{for each separating node $\sta e$},
\end{align*}
where, for each separating node, $d(\sta e)$
stands for  the degree of $\pi^*F$ on one of the connected components
of the partial normalisation  of $\sta C$ at $\sta e$.\qed
\end{thm}
\subsection{An exact sequence relating $\Pic_{\stC}[r]$ and $\Pic_{\coa{\stC}}[r]$.}\label{sect:coast}
In fact if $\pi\colon\sta C\to \coa{\stC}$ is
a twisted curve over $k$,
for which \emph{all
nodes have stabilisers of order $r$}, once an orientation
of the corresponding dual graph $\Gamma$ is chosen,
we have the following exact sequence from \cite[3.0.8]{Ch_mod}
\begin{equation}\label{eq:seq}
1\to \Pic_{\coa{\stC}}[r]
\xrightarrow{\pi^*}
\Pic_{\stC}[r]
\xrightarrow{\sta j^*} C_1(\Gamma,\ZZ/r\ZZ)
\xrightarrow {\partial}
C_0(\Gamma,\ZZ/r\ZZ)\xrightarrow {\varepsilon}\ZZ/r\ZZ\to 1
\end{equation}
where the homomorphisms are defined as follows. The differential $\partial$ is
the boundary homomorphism with respect to the  orientation of
$\Gamma$, and $\varepsilon$ denotes the augmentation homomorphism
$(h_v)_V\mapsto \sum_V h_v\in \ZZ/r\ZZ$. Finally
$\sta j^*$ can be regarded as the pullback to the singular locus via
$\sta j\colon {\rm Sing}_{\stC} \to \stC$ or, more explicitly,
as the homomorphism mapping $\sta L\in \Pic_{\stC}[r]$ to
$({q_e})_{e\in E}\in C^1(\Gamma,\ZZ/r\ZZ)$, where $q_e\in \{0, \dots, r-1\}$ is defined as in
\S\eqref{sect:monod}.
\subsection{The automorphism group of a twisted curve.}
\label{rem:lpghost}
Let $\pi\colon \stC\to \coa{\sta C}$ be a twisted curve over $k$ with
all stabilisers of order $r$.
In \cite[Prop.~7.1.1]{ACV}
the group
$\Aut(\stC,\coa{\sta C})$ of automorphisms
of $\stC$ fixing the coarse space $\coa{\sta C}$
is explicitly calculated:
for any  twisted curve $\stC$ over $k$, we have
an isomorphism
$$\Aut(\stC,\coa{\sta C})\cong ( \pmmu_{r})^E.$$

In fact \cite[Thm.~7.1.1]{ACV}
shows that we can choose a set of
independent generators as follows. For $\sta e\in \stC$, we have an automorphism
$\sta g\in \Aut(\stC,\coa{\sta C})$
such that the restriction of
$\sta g$ to ${\stC\setminus \{\sta e\}}$ is the identity,
and the local picture at $\sta e$ is given by
\begin{equation}\label{eq:gloc}
(z,w) \mapsto (z,\xi_r w).
\end{equation} operating on the scheme $\Spec (k[z,w]/(zw))$.
(Note that morphisms of stacks are given up to natural transformations. In this way no
branch has been privileged. The
$1$-automorphism in the local picture above is in fact, $2$-isomorphic to
the morphism $(z,w)\mapsto (\xi^{h}_rz,\xi_r^{1-h}y)$
for any $h=0,\dots, r-1$.)

\subsection{The action of the automorphism group on the Picard group.}
We describe the action on $\Pic_{\stC}[r]$ by pullback via the above
automorphism $\sta g$.
We choose an orientation of the dual graph of $\stC$.
We choose a primitive $r$th root of unity $\xi_r$, and in
\eqref{eq:seq} we identify $\pmmu_r$ with $\ZZ/r\ZZ$ accordingly.
We assume that, with this choice of $\xi_r$, the
local picture of $\sta g$ at $\sta e$ is given by
$(z,w) \mapsto (z,\xi_r w)$ as in \eqref{eq:gloc}.
We denote by $\stC^\nu$ the normalisation
$\nu\colon \stC^\nu\to \stC$
and recall that we have the exact sequence
$1\to \pmmu_r\to C_0(\Gamma,\pmmu_r) \xrightarrow{\ \delta \ } C_1(\Gamma,\pmmu_r)
\to \Pic_{\stC}[r]\xrightarrow{\ \nu^* \ } \Pic_{\stC^\nu}[r].$

\begin{pro}\label{pro:aut_lb}
The pullback via the automorphism $\sta g
\in \Aut(\stC,\coa{\stC})$
satisfying \eqref{eq:gloc}
can be written as
\begin{align*}
\label{eq:adjust}
\sta g^*\colon \Pic _{\stC}[r] &\to \Pic _{\stC}[r]\\
\stL&\mapsto  \stL\otimes \sta T_{\sta L},
\end{align*}
where $\sta L\mapsto \sta T_\sta L$ is the composite homomorphism
of $\Pic_{\stC}[r]\to C_1(\Gamma,\pmmu_r)$ and
$C_1(\Gamma,\pmmu_r)\to \Pic_{\stC}[r]$ fitting in the exact sequences
recalled above and in \eqref{eq:seq}.
\end{pro}
\begin{proof}
The claim follows from \cite[Prop.~2.5.4]{Ch_mod} which identifies the line bundle $\sta T_\stL$
with the composite morphism described above. Indeed
$\sta T_\stL$ is the sheaf of regular functions $f$ on $\stC^\nu$ whose values on the preimages
$\sta p_+$ and $\sta p_-$ of a node $\sta e$
are related by $f(\sta p_+)=\xi_r^{q_e}f(\sta p_-)$. This is precisely
the image via  $C_1(\Gamma,\pmmu_r)\to \Pic_{\stC}[r]$ of
$(\xi_r^{q_\sta e})_{\sta e\in E}\in C_1(\Gamma,\pmmu_r)$, with $q_{\sta e}$ defined as in \S\eqref{sect:monod}
for each node $\sta e$.
\end{proof}
\begin{rem}\label{rem:torsor}
In fact $\sta g$ acts by pullback on the entire Picard group. By \cite[Prop.~2.5.4]{Ch_mod}
the action on $\Pic_{\stC}$ is again the tensor product
of the identity and a morphism $\stL\mapsto \sta T_\stL$ given by
composing
$\sta j^*\colon \Pic_{\stC}\to C_1(\Gamma,\pmmu_r)$ with
$C_1(\Gamma,\pmmu_r)\to \Pic_{\stC}$.\end{rem}

\section{N\'eron $d$-models}\label{sect:dmodels}
\subsection{The definition of N\'eron model.}
For any
$K$-scheme $X_K$, a scheme $Y$ over $S$
is an \linebreak \emph{$S$-model of $X_K$} if
 the generic fibre is $X_K$.
There is an abundance of $S$-models; on the other hand,
the N\'eron model is a canonical $S$-model
smooth, separated,
and of finite type, and
satisfying  a universal
property, the N\'eron mapping property,
which determines it uniquely, up to a canonical isomorphism:\\
For each smooth $S$-scheme
$Y\to S$ and each $K$-morphism
$u_K\colon Y_K\to X_K$, there is a unique
$S$-morphism $u\colon Y\to X$
extending $u_K$.\\
The N\'eron model commutes with \'etale base changes,
with the passage
 to henselisation \cite{Re_Pic}.
\subsection{Existence of N\'eron models.}\label{sect:exist}
The existence of the
N\'eron models of group schemes and torsors is proven
in  \cite[4.3/Thm.~6]{BLR},
\cite[4.4/Cor.~4]{BLR}
under a boundedness assumption \cite[1.1/Def.~2]{BLR}.
In this paper we only consider N\'eron models of proper
$K$-group schemes and proper $K$-torsors, so that
the boundedness assumption is automatically satisfied.
We have the following statements.

Let $G_K\to \Spec K$ be a group scheme.
Assume that it is smooth, of finite type, and proper over
$K$. Then, a N\'eron model
$G$ on $S=\Spec R$ exists and is unique, up to a canonical isomorphism.
The structure of $G_K$
as a group $K$-scheme extends uniquely to a structure of $G$
as a group $S$-scheme.

Furthermore, assume that
$T_K$ is also smooth, of finite type, and proper over $K$ and
is a torsor on $\Spec K$ under $G_K$.
Then, a N\'eron model $T$ on $S=\Spec R$ exists
and is unique, up to a canonical isomorphism.
If $T\to S$ is surjective, the structure of $T_K$
as a torsor on $\Spec K$ under $G_K$ extends
uniquely to a structure of $T$
as a torsor on $S$ under $G$, see \cite[6.5/Cor.~3.4]{BLR}.

\subsection{The definition of N\'eron $d$-model.}
We place the notion of N\'eron model on a
stack-theoretic base. Instead of $S$,
for a suitable positive
integer $d$ prime to $\cha (k)$,
we take as a base the quotient stack
$S[d]=[\Spec R_d/\pmmu_d]$, where
$R_d$ is equal to $R[\wt \pi]/(\wt\pi^d-\pi)$
for a uniformiser $\pi$ of $R$
and $t\in \pmmu_d$ acts on $R_d$
by $t(\wt \pi)=t\wt \pi$ and fixes $R$. In this way, we have
$$\Spec K\xrightarrow{\ i\ } S[d]\xrightarrow{\ p \ } \Spec R,$$
where $i$ is an open and dense immersion and
$p$ is the (proper) morphism to the coarse space.
\begin{defn}[{$S[d]$}-model]
For any $K$-scheme $X_K$,
a \emph{representable} morphism of stacks\linebreak
$\sta Y\to S[d]$
is an \emph{$S[d]$-model of $X_K$}
if its generic fibre $\sta Y\times _{S[d]}\Spec K$ is $X_K$.
\end{defn}
We recall that, since we are
working inside the $2$-category of algebraic stacks, an \linebreak
$S[d]$-morphism from $\sta f\colon \sta X\to S[d]$
to
$\sta f'\colon \sta X'\to  S[d]$
is a morphism
$\sta g\colon \sta X\to \sta X'$
alongside with a $2$-isomorphism
$\sta g\circ \sta f' \Rightarrow \sta f$.
Note, however, that the $2$-isomorphism is
uniquely determined, because
\emph{any automorphism of a  representable
smooth morphism $\sta f\colon \sta X\to S[d]$
is trivial} (Lemma 4.2.3 of  \cite{AV} applies
since $\sta f$ maps the open dense
subscheme $\sta X_K$ of $\sta X$
into the open dense subscheme $\Spec K$ of $S[d]$).

\begin{defn}[N\'eron $d$-model]\label{defn:dner}
Let $\sta X_K$ be a smooth and separated $K$-scheme of finite type.
For $d$ prime to  $\cha(k)$, a
N\'eron $d$-model is an $S[d]$-model $\sta X$ of $X_K$,
which is smooth, separated, and of finite type,
and which satisfies the following universal property
analogue to the N\'eron mapping property:

For each representable and smooth
morphism of stacks $\sta Y\to S[d]$, and each $K$-morphism
$\sta u_K\colon \sta Y_K\to \sta X_K$, there is an
$S[d]$-morphism $\sta u\colon \sta Y\to \sta X$,
which extends $\sta u_K$ and is unique, up to
a unique natural transformation
\[\xymatrix@R=0.5cm{
\sta Y_K\ar[d]_{\sta{u}_K}\ar[r]& \sta Y\ar[d]_{\sta{u}}\\
\sta X_K\ar[d]      \ar[r]& \sta X\ar[d]\\
\Spec K\ar[r]        & S[d].
}\]
\end{defn}

\begin{rem}\label{rem:unicity}
It follows from the Definition \ref{defn:dner}, that any two
N\'eron $d$-models of $X_K$ are isomorphic and the isomorphism
is unique up to a unique natural transformation.
\end{rem}

\begin{pro}\label{pro:relner}
Let $d$ be invertible in the residue field.
For $\pi$ a uniformiser of $R$,
write $\wt R$ for $R[\wt\pi]/(\wt\pi^d-\pi)$,
 $\wt K$ for the corresponding field of fractions,
and $\wt S$ for $\Spec \wt R$, with the natural $\pmmu_d$-action.
Let $X_K$ be a
smooth and separated $K$-scheme of finite type,
and let $X_{\wt K}=X_K\otimes _K \wt K$ be the corresponding $\wt K$-scheme.
\begin{enumerate}
\item[I.] If $\sta N$ is the N\'eron $d$-model of $X_K$,
then the $\wt S$-scheme $\wt {N}$ fitting in
\[
\xymatrix@R=0.2cm{
\wt {N}\ar[rr]\ar[dd]&&\sta N\ar[dd]\\
&\square&\\
\wt S\ar[rr]&&S[d]
}
\]
is the
N\'eron model of $X_{\wt K}$.
\item[II.] Conversely, assume that the scheme
$X_{\wt K}$ has a
N\'eron model $N(X_{\wt K})$ on $\wt S$. Then,
there is a natural $\pmmu_d$-action on $N(X_{\wt K})$ together with  a
$\pmmu_d$-equivariant morphism $N(X_{\wt K})\to \wt S$, and
the corresponding morphism of stacks
$$[N(X_{\wt K})/\pmmu_d]\to S[d]$$ is
the N\'eron $d$-model of $X_K$.
\end{enumerate}
\end{pro}
\begin{proof}
Since $\wt S\to S[d]$ is \'etale, point (I) is a straightforward consequence of
the compatibility of
the formation of N\'eron models with finite \'etale base change, \cite{Re_Pic}.

For (II), note that
$\pmmu_d$ acts on $N(X_{\wt K})$, because it acts on
$X_{\wt K}$ and, by the N\'eron mapping property,
the action extends to $N(X_{\wt K})$.
Therefore, $[N(X_{\wt K})/\pmmu_d]$ is an $S[d]$-model
of $X_{K}=[X_{\wt K}/\pmmu_d]$.

In order to check the universal property of Definition \ref{defn:dner},
consider a smooth and representable
morphism $\sta Y\to S[d]$
and a morphism $\sta u_K\colon \sta Y_K\to X_K$.
Note that $\sta Y\to S[d]$ can
be regarded as a $\pmmu_d$-equivariant smooth $\wt S$-scheme $\wt Y$:
indeed, $\wt Y$ is defined as $\sta Y\times _{S[d]} \wt S$,
which is a scheme by the representability assumption, and
the $\pmmu_d$-action is defined by pullback of
$\pmmu_d\times \wt S\to \wt S$.
In this way, $\sta u_K$ lifts to a $\pmmu_d$-equivariant
morphism $\wt Y_{\wt K}\to X_{\wt K}$.
By the N\'eron mapping property for $N(X_{\wt K})$,
we have a $\pmmu_d$-equivariant  $\wt S$-morphism
$\wt u\colon \wt Y\to N(X_{\wt K})$
extending $\wt Y_{\wt K}\to X_{\wt K}$.
So,
$\sta u\colon \sta Y=[\wt Y/\pmmu_d]\to [N(X_{\wt K})/\pmmu_d]$ extends
$\sta u_K\colon \sta Y_K\to X_K$.

Finally, take a morphism
$\sta u'\colon \sta Y \to [N(X_{\wt K})/\pmmu_d]$
such that $\sta u'\otimes K$
coincides with $\sta u_{K}$ on $\sta Y_K$. In fact,
$\sta u'\otimes K$  and  $\sta u_{K}$
are lifted by morphisms from $\wt Y_{\wt K}$ to $X_{\wt K}$
and the two liftings coincide after composition with the action of an
element of $\pmmu_d$.
By the N\'eron mapping property, the extension to $\wt Y$ also coincide
up to the action of an element of $\pmmu_d$; this means that
the morphisms of stacks $\sta u$ and $\sta u'$ are isomorphic
up to a unique natural transformation.
\end{proof}

\begin{rem}\label{rem:existence}
By Proposition \ref{pro:relner}, the existence of N\'eron $d$-models is
guaranteed under the properness assumptions of
\S\eqref{sect:exist}.
Furthermore, the N\'eron $d$-model of
a group $K$-scheme is equipped with a unique structure of
group stack on $S[d]$, and the N\'eron $d$-model of
a $K$-torsor is equipped with a unique structure of
torsor on $S[d]$ if it surjects on $S[d]$.
\end{rem}
\begin{notn}
Whenever $d$ is invertible in the residue field $k$ and
$X_K$ satisfies the hypothesis of \S\eqref{sect:exist},
the N\'eron $d$-model of $X_K$ exists and we denote it by $$N_d(X_K)\to S[d].$$
\end{notn}

\section{Finite N\'eron $d$-models of the $r$-torsion of the Picard group}\label{sect:group}
\subsection{The finiteness of the N\'eron $d$-model with respect to $d$.}
In the following theorem we identify the finite N\'eron $d$-models and those
who are represented by $r$-torsion line bundles on a twisted semistable reduction.
\begin{thm}\label{thm:group}
Let $r>2$ be an integer prime to $\cha(k)$, and
let $C_K$ be a smooth curve of genus $g\ge 2$.
Assume that the group $K$-scheme $G_K=\Pic_{C_K}[r]$ is
tamely ramified on $K$, i.e. all its points correspond to
tame extensions of $K$.

Then, there exist three integer and positive invariants
$m_1$, $m_2$, and $m_3$
of $C_K$ satisfying
\begin{align*}
&m_2\in m_1\ZZ,\\
&            m_3\in m_2\ZZ,\\
& rm_1\in m_3\ZZ,
\end{align*}
and the following conditions.
\begin{enumerate}
\item
There is a semistable reduction of $C_K$ on $S[d]$ if and only if
$d$ is a multiple of $m_1$.
\item
The N\'eron $d$-model $N_d(G_K)$ is a finite group scheme
if and only if $d$ is a multiple of $m_2$.
\item
The N\'eron $d$-model $N_d(G_K)$ is a finite group scheme
and represents the $r$th roots of $\Ocal$ on a twisted
semistable reduction $\stC$ of $C_K$ on $S[d]$
if and only if $d$ is a multiple of $m_3$.
\end{enumerate}
If $C_K$ admits a semistable reduction over $R$, then
the condition $r>2$ is superfluous and the condition that
$G_K$ is tamely ramified is always true. In this case, for any $r$ prime to $\cha(k)$, we have $m_1=1$,
$m_2\mid m_3$, and $m_3\mid r$.
\end{thm}

\begin{rem}\label{rem:credits}
Point (1) follows easily from
a version of the theorem of semistable reduction (the argument is given in
\cite[\S5]{De} by Deschamps and is attributed to Raynaud).
We also point out that $m_2$ is a multiple of $m_1$ is merely a reformulation
of a criterion due to Serre
(see \cite[\S5]{De}). These facts are reviewed in the course of the proof
because they are needed in the rest of the paper. \end{rem}

\begin{proof}[Proof of Theorem \ref{thm:group}]
We prove (2), from which we deduce (1). Finally, we show point (3).

Point (2) affirms the existence of $m_2$ such
that $N_d(\Pic_{C_K}[r])$ is finite if and only if $m_2$ divides $d$.
This is a simple fact in Galois theory,
which we need later; we recall it in the following lemma.
\begin{lem}\label{lem:finKgr}
Let $G_K$ be a tamely ramified finite $K$-scheme. Then
there is an integer and positive invariant $m(G_K)$
such that the N\'eron $d$-model of $G_K$ is finite if and only if
$d$ is a multiple of $m(G_K)$.
\end{lem}
\begin{proof}
Denote by $\ol K$
a separable algebraic closure
of $K$; we write $\ol G$
for $G_K\otimes \ol K$.
By Proposition \ref{pro:relner}, we only need to
determine the integers $d$ for which
$R_d=R[\pi']/(\pi'^d-\pi)$
satisfies the following condition:
the N\'eron model of
the pullback $G_d$ of $G_K$ on
the corresponding valuation field $K_d\supset R_d$
is finite on $R_d$.
This happens if and only if $G_d$
is not ramified.

By descent theory, there is a natural morphism
\begin{equation}\label{eq:gal}
\Gal(\ol K/K)\to \Aut(\ol G).
\end{equation}
Note that  the $K_d$-scheme $G_d$
is not ramified
if and only if
$\Gal(\ol K/K_d)$ is contained in the kernel of
the above morphism.
Since $G_K$ is tamely ramified, the image
of \eqref{eq:gal} is a finite cyclic group
whose order is prime to $\cha(k)$.
Let $m(G_K)$ be such order.
Then, $G_d$ is not ramified on $R_d$
if and only if $d$ is a multiple of $m(G_K)$.
\end{proof}

Point (1) follows from the following statement and from the application of Lemma \ref{lem:finKgr}.
\begin{lem}[{\cite[\S5]{De}}]\label{lem:smallgroup}
There exists a group $K$-scheme $E_K$ associated to $C_K$ and to $r>2$
satisfying the following property.
The curve $C_K$ has a
semistable reduction  on $S[d]$
if and only if the N\'eron $d$-model of $E_K$ is finite.
\end{lem}
\begin{proof}
In order to define $E_K$, one needs to introduce
a finite Galois extension $K\subset K'$
for which $C_K\otimes  K'$ has semistable reduction over $R'$,
the integral closure of $R$ in $K'$ (such an extension exists by the
theorem of semistable reduction).
Then, there is a unique group $E_K$ satisfying the following conditions.

Let $U'$ be the group $R'$-scheme representing the functor
of line bundles of degree zero \emph{on all irreducible
components} of the fibres of the semistable minimal regular model of $C_K$ on $\wt R$.
Let P be the $r$-torsion
subgroup of $U'$.
The group $R'$-scheme $P$
is the disjoint union $P=P^{\rm gen}\sqcup P^{\rm fin}$, where
$P^{\rm gen}$  is a component contained in the generic fibre and a
$P^{\rm fin}$ is a finite group scheme over the base ring $R'$.
Then, $E_K$ is determined by descent: the group $K$-scheme $E_K$ is the
group satisfying $E_K\otimes K'=P^{\rm fin}\otimes K'$.
The fact that $P^{\rm fin}\otimes K'$ descends and that this definition
does not depend on the Galois extension $K'$
is shown in  {\cite[\S5]{De}}.
It remains to show that the equivalence holds.

By construction of $E_K$, if
$C_K$ has stable reduction on $S[d]$,
then $N_d(E_K)$ is finite. Indeed,
since the definition of $E_K$ does not depend on the extension, we can
define $E_K$ by descent from $K_d=K[\pi']/(\pi'^d-\pi)$ to $K$. The group
$K_d$-scheme
$E_K\otimes K_d$ is the generic fibre of a finite
group scheme $P^{\rm fin}$ over
 the integral closure $R_d$ of $R$ in  $K_d$.
 Then, the N\'eron model of  $E_K\otimes K_d$ over $R_d$ is finite.
 By Proposition \ref{pro:relner},
this is equivalent to say that
the N\'eron $d$-model of $E_K$ is
finite.

Conversely we prove that, if the N\'eron model of $E_K\otimes K_d$ is finite over $R_d$,
then the minimal regular model of
$C_K\otimes K_d$ is semistable over $R_d$.
Recall that we can choose a
finite Galois extension $K'$ containing $K_d$ such that
the minimal regular model of $C_K\otimes K'$
over the integral closure $R'$ of $R$ in $K'$ is semistable.
The claim is equivalent to showing that the jacobian variety of $C_K\otimes K_d$ has
semistable reduction over $R_d$. This is the case
if there is a semistable reduction of
the jacobian variety of $C_K\otimes K'$
which descends from $R'$ to $R_d$.

Now, note that such a semistable reduction is provided by the
group $R'$-scheme $U'$ defined above.
Indeed, it descends from $R'$ to $R_d$ because the action of
$\Gal(K'/K_d)$ on its special fibre is trivial. In order to see that, via Serre's lemma \cite{Serre},
it is enough to show that $\Gal(K'/K_d)$ acts trivially on the special fibre of the
kernel of the multiplication by $r$: the scheme $P$ introduced above.
In fact, the special fibre of $P$ coincides
with the special fibre of $P^{\rm fin}$ and the action of $\Gal(K'/K_d)$
on the special fibre of $P^{\rm fin}$ is trivial as soon as $P^{\rm fin}$ descends
to $R_d$. This is indeed a consequence of the fact that
the N\'eron model of $E_K\otimes {K_d}$ is finite over $R_d$.
\end{proof}
Note that, since $E_K$ is
a subgroup of $\Pic_{C_K}[r]$ by definition, the above statement
also shows that
$m_1$ divides $m_2$.

Point (3) determines all
indices $d$ satisfying the property P($d$):
the N\'eron $d$-model $N_d(G_K)$ is a finite group scheme
and represents the $r$th roots of $\Ocal$ on a twisted
semistable reduction $\stC$ of $C_K$ on $S[d]$.
This is a consequence of Theorem \ref{thm:cond},
which states that a twisted curve over $k$ of genus $g$ has $r^{2g}$ $r$-torsion
line bundles if
and only if the order of the stabiliser of each nonseparating node
is a multiple of $r$. Applying this fact, we can immediately show that
for any $d\in rm_1\ZZ$ the condition P$(d)$ is satisfied.
Indeed, let $C_{m_1}$ be the stable reduction of $C_K$ over $S[m_1]$.
Then, the pullback to $S[d]$ is a stable reduction
$C_{d}$ of $C_K$. Note that all nodes have thicknesses in $r\ZZ$ because
the base change of  $\{zw=s^l\} $ via $s\mapsto s^h$
yields $\{zw=s^{hl}\}$ for all positive integers $h$ and $l$.
Then, by Remark \ref{rem:modifsst}, on $S[d]$ there is a twisted curve $\stC_{d}$
whose nonseparating nodes
have stabilisers whose order lies in $r\ZZ$.
Finally, by
Theorem \ref{thm:cond} the group scheme formed by $r$th roots of
$\Ocal$ on $\wt \stC$ is finite and
is the N\'eron $d$-model of its
generic fibre.

This shows that $rm_1\ZZ$ is included in the set of
indices for which P($d$) is satisfied.
Assume that we modify the statement P$(d)$: we
require   that P$(d)$ is true and
also that the coarse space $\coa{\sta C}$ of
$\stC$ is stable. Let us call this property P$(d)$'. Then
the argument given above shows without change that
the indices $d$ satisfy ${\rm P}(d)$' if and only if they are
multiple of $m_1 r/\gcd\{r,t\}$,
where $t$ is
the greatest common divisor of
the thicknesses of nonseparating nodes of the stable model $C_{m_1}$ over $S[m_1]$.
In fact, this proves point (3) entirely, because we have
${\rm P}(d)\Leftrightarrow {\rm P}(d)'$. Indeed,  ${\rm P}(d)'\Rightarrow {\rm P}(d)$ is obvious.
Conversely, if
P$(d)$ is true, then there is a twisted semistable reduction ${\stC}$ of $C_K$ on $S[d]$
such that $r$ divides the height of
all nonseparating nodes in the coarse space $\coa{\stC}$.
After contraction of all projective lines
with only two nodes in the special fibre of $\coa{\stC}$, we
obtain the stable model for which
the thicknesses of the nonseparating nodes still belong to  $r\ZZ$.
This elementary fact easily follows from the fact
that, in the semistable minimal regular model obtained by desingularising
$\coa{\stC}$,
the number of nonseparating nodes in each chain of $-2$-curves is a multiple of $r$.
\end{proof}

The proof of Theorem \ref{thm:group}
provides some characterisations of
the invariants $m_1$, $m_2$, and $m_3$, which we
state explicitly.

\subsection{The indices $m_1$ and $m_2$.}
By Lemma \ref{lem:finKgr}, we have
$$m_1=\size{\im(d_{E_K})},$$
where $E_K$ is the
subgroup of $G_K$ satisfying Lemma \ref{lem:smallgroup} and
$d_{E_K}$ is the morphism \linebreak $\Gal(\ol K/ K)\to \Aut (E_K\otimes \ol K)$
for a separable closure
$\ol K$ of $K$.
Similarly, we have
$$m_2=\size{\im(d_{G_K})}.$$

\subsection{The ratio $m_3/m_1$.}\label{sect:m3m1}
The ratio $m_3/m_1$
can be expressed in terms of the geometry of the stable reduction
of $C_K$ on $S[m_1]$.
Indeed, we have
\begin{equation}\label{eq:m1m3}
m_3/m_1= r/\gcd\{r, t\},\end{equation}
where $t$ is the greatest common divisor of the thicknesses of the
nonseparating nodes appearing in the stable reduction of
$C_K$ on $S[m_1]$.
As a consequence we have the following statement.
\begin{cor}\label{cor:repr}
The N\'eron model
of $\Pic_{C_K}[r]$ over $R$ is finite
and represents the functor of
$r$-torsion line bundles on a twisted semistable
reduction of $C_K$ over $R$ if and only if $C_K$ admits a stable model
in which the thickness of each nonseparating nodes in the
special fibre is a multiple of $r$.
\end{cor}
\subsection{The ratio $m_2/m_1$.}\label{sect:m2m1}
The argument of the proof of Theorem \ref{thm:group} allows us
to express the ratio $m_2/m_1$ in a similar way.
\begin{pro}\label{pro:m2m1}
We have
\begin{equation}\label{eq:m1m2}
m_2/m_1 = r/\gcd\{r, c\},\end{equation}
where $c$ is the greatest common divisor of the  number  of
edges common to two  circuits in the dual graph
of the special fibre of the semistable minimal regular model of
$C_K$ on $S[m_1]$.\end{pro}
\begin{proof}
We already showed the relations
$m_1\ZZ\subseteq m_2\ZZ\subseteq m_3\ZZ\subseteq rm_1\ZZ$.
It remains to show that the divisor $q$ of $r$ satisfying
$rm_1=qm_2$ equals $\gcd\{r, c\}$.
The claim can be equivalently restated as follows. Assume $m_1=1$; then
$\gcd\{r, c\}$ is the highest among the
divisors $q$ of $r$ for which
the finite N\'eron $r$-model descends to $S[r/q]$. In other terms, $\gcd\{r, c\}$ is
the highest among the divisors $q$
of $r$ for which $\pmmu_{q}\subseteq \pmmu_{r}$ acts trivially on the special fibre
of $N_{r}(G_K)$.

We can realise $N_{r}(G_K)$ as the group of $r$-torsion line bundles on
a twisted semistable reduction of $C_K$ on $S[r]$
whose fibers only contain stabilisers of order $r$ on each node.
We construct this twisted semistable reduction as we did in the
proof of Theorem \ref{thm:group}:
by iterating the construction illustrated in Remark \ref{rem:modifsst}.
In this way, we can assume that we have a
twisted semistable reduction of $C_K$ on $S[r]$
such that the coarse space descends on $S$ and is the semistable
minimal regular model of $C_K$.
We regard this twisted curve on $S[r]$ as a \linebreak $\pmmu_{r}$-equivariant
twisted curve $\wt \stC$ over the discrete valuation ring $\wt R=R[\wt \pi]/(\wt\pi^r-\pi)$.

The action of $\pmmu_r$ on $\wt\stC$
is trivial on the special fibre of the coarse space
$\coa{\wt\stC}$, because
the coarse space descends to $S$.
In this way, the action of $\pmmu_r$  on $\wt \stC$ is characterised by
its local picture at the nodes.
We recall that the local picture of $\wt \stC$ at a node is
$[U/\pmmu_r]$, where
$U$ is $\{zw=\wt \pi\}$, and the quotient stack
is defined by the $\pmmu_r$-action $(z,w,\pi)\mapsto (\xi_rz,\xi_r^{-1} w, \wt\pi)$.
Therefore, we can write
$\{xy=\wt \pi^r\}$ for the local picture of $\coa{\stC}$
(with $x=z^r$ and $y=z^r$),
and we notice that $\{xy=\wt \pi^r\}$ is the pullback via $\wt \pi\mapsto \wt\pi^r$
of $\{xy=\pi\}$, the local picture of the minimal regular model.
Now, $\xi_r\in \pmmu_r$ operates nontrivially on $[U/\pmmu_r]$
but trivially on the special fibre of the coarse scheme. Therefore,
$\xi_r\in \pmmu_r$ fixes $zw-\wt \pi$, $x=z^r$, $y=q^r$, and acts by multiplication of $\wt \pi$.
This implies that $\xi_r$ operates on $[U/\pmmu_r]$ as
$$(z,w,\pi)\mapsto (\xi_r^a z,\xi_r^b w, \xi_r\wt\pi)\quad \quad \text{with \ \ } a+b\equiv 1\mod r.$$
Note that, up to natural transformation\footnote{The natural
transformations have
 equation $(z,w,\pi)\mapsto (\xi^h_rz,\xi_r^{-h} w, \wt\pi)$},
this local picture identifies a single morphism $\sta g$ whose local picture at all nodes is given by
\begin{equation}\label{eq:g}
(z,w,\pi)\mapsto (z,\xi_r w, \xi_r\wt\pi).\end{equation}
In this way, the action of
$\pmmu_r$ on the special fibre $\wt\stC_k$
is generated by the automorphism $\sta g$ studied in Proposition \ref{pro:aut_lb}.

Using Proposition \ref{pro:aut_lb}, we conclude that
the highest divisor $q$ of $r$ for which $\pmmu_{q}$ in $\pmmu_r$
acts trivially on $\Pic_{\wt\stC_k}[r]$ is the highest index $q$ for which
the endomorphism of $\Pic_{\wt\stC_k}[r]$ given by
$\stL \mapsto \sta T_{\stL}^{\otimes r/q}$ vanishes.
Proposition \ref{pro:aut_lb} claims that
$\stL \mapsto \sta T_{\stL}$ fits in the following diagram where the vertical and the horizontal
sequences are exact.
\begin{equation}\label{eq:Xdiag}
\xymatrix{&\Pic_{\stC}[r]\\
\Pic_{\stC}[r]\ar[ur]^{\stL \mapsto \sta T_{\stL}}\ar[r]^{\sta j^*}& C_1(\Gamma,\pmmu_r)\ar[r]^{\partial_r}\ar[u]&C_0(\Gamma,\pmmu_r)\\
& C_0(\Gamma,\pmmu_r)\ar[u]_{\delta_r}.}\end{equation}
Therefore, $\stL \mapsto \sta T_{\stL}$ vanishes if and only if
the dual graph
$\Gamma$ of $\wt\stC_k$ satisfies $\ker(\partial_r) \subseteq \im (\delta_r)$.
More generally, it is easy to see that $\stL \mapsto \sta T_{\stL}^{\otimes r/q}$ vanishes if and only if
$\Gamma$ satisfies $\ker(\partial_q) \subseteq \im (\delta_q)$.
This completes the proof, because imposing $\ker(\partial_q) \subseteq \im (\delta_q)$
is equivalent to imposing the condition
that the
number  of edges common to any two  circuits is always a multiple of $q$. See Section \ref{sect:circuits}
for a proof.
\end{proof}

\subsection{The finiteness criterion for $N(\Pic_{C_K}[r])$.}
As a consequence of the previous results, we can determine whether the
N\'eron model of
$\Pic_{C_K}[r]$ is finite over $R$
simply by looking at the special fibre of the
minimal regular model of $C_K$. .
\begin{cor}\label{cor:group}
Let $r>2$ be prime to the residual characteristic, and let $C_K$
be a smooth curve of genus $g\ge 2$.

The N\'eron model of
$\Pic_{C_K}[r]$ is finite if and only if
$C_K$ admits a semistable minimal regular model for which
the dual graph of the special fibre satisfies the following property:
the (signed) number of edges common to any two  circuits is always a multiple of $r$.
\end{cor}
\begin{proof}
The N\'eron model $N(\Pic_{C_K}[r])$ is finite if and only if
$m_2=1$.
By \eqref{eq:m1m2}, $m_2=1$ if and only if
$m_1=1$ and $c\in r\ZZ$.
This is precisely the statement.
\end{proof}
\begin{rem}
Note that the two conditions related in the statement of Corollary \ref{cor:group} hold only if
the group $K$-scheme
$\Pic_{C_K}[r]$ is tamely ramified.
Therefore, do not need to  assume this condition in the statement.
\end{rem}
\begin{rem}\label{rem:r=2} The criterion given above holds with the same proof
for $r=2$ if we make the assumption that
a semistable reduction of $C_K$ over $R$ exists. Otherwise, in the statement, ``if and only if''
should be replaced by ``if''.
\end{rem}
\subsection{Some examples.}\label{sect:examples}
We compare the values of  $m_1$, $m_2$, and $m_3$
in some concrete examples.
We always consider a \emph{stable} curve $C$ over $R$ satisfying the following conditions:
\begin{enumerate}
\item the generic fibre is a smooth curve $C_K$ over $K$,
\item $C$ is regular (in other words, $C$ is the minimal regular model of $C_K$).
\end{enumerate}
The dual graph of the special fibre is denoted by $\Gamma$.
We list the values of $m_1,m_2$, and $m_3$
in the following table, where we assume
for simplicity that $r$ is a positive multiple of $4$.
Note that $m_1$ always equals 1 because
the curve $C_K$ admits a semistable reduction. Furthermore, $m_3$ equals
$r$ because all nodes of the stable model have thickness $1$.
On the other hand $m_2$ is determined by the (signed) numbers
of edges shared by two
circuits in the graph. The greatest common divisor of these numbers is
$1,2,4,2,4$, and $2$ in the four cases considered; so $m_2$ equals $r,r/2,r/4,r/2,r/4$, and $r/2$, respectively.

\begin{center}
\begin{tabular}{|c|c c c c|}

\hline
dual graph $\Gamma$ of the special fibre of $C$
&
$m_1$
&
$m_2$
&
$m_3$
&
\\
\hline
\xymatrix@R=.5pc{
*{ \,\bullet} \ar@{-}@(ur,dr) & {}\\
{}
}

&
\xymatrix@=-.5pc{
\\
1}

&
\xymatrix@=-.4pc{
\\
{r}}
&
\xymatrix@=-.4pc{
\\
r}
&
\\
\hline

\xymatrix{
*{ \,\bullet} \ar@{-}@/_/[rd] \\
& *{ \,\bullet} \ar@{-}@/_/[lu]
}

&
\xymatrix@=.3pc{
\\
1}

&
\xymatrix@=.18pc{
\\
{r}/{2}}
&
\xymatrix@=.45pc{
\\
r}
&
\\

\hline

\xymatrix{
*{ \,\bullet} \ar@{-}[r] &*{ \,\bullet} \ar@{-}[d] \\
*{ \,\bullet} \ar@{-}[u] & *{ \,\bullet} \ar@{-}[l]
}

&
\xymatrix@=.3pc{
\\
1}

&
\xymatrix@=.18pc{
\\
{r}/{4}}
&
\xymatrix@=.45pc{
\\
r}
&
\\

\hline
\xymatrix@=1.6pc{
&*{\bullet} \ar@{-}[dl]\ar@{-}[d] \ar@{-}[d]\ar@{-}[dr] \\
*{\bullet} \ar@{-}[dr] &*{\bullet} \ar@{-}[d] &*{\bullet} \ar@{-}[dl] \\
&*{\bullet} \\}

&
\xymatrix@=1.1pc{
\\
1}

&
\xymatrix@=0.98pc{
\\
{r/2}}
&
\xymatrix@=1.25pc{
\\
r}
&
\\
\hline

\xymatrix@=1pc{
&*{\bullet} \ar@{-}[dl] \ar@{-}[dr]&&&*{\bullet} \ar@{-}[dl] \ar@{-}[dr] \\
*{\bullet} \ar@{-}[dr] &    &*{\bullet} \ar@{-}[dl]\ar@{-}[r] &*{\bullet} \ar@{-}[dr] &     &*{\bullet} \ar@{-}[dl]\\
&*{\bullet}  &&&*{\bullet} \\}

&
\xymatrix@=.4pc{
\\
1}

&
\xymatrix@=0.28pc{
\\
{r}/{4}}
&
\xymatrix@=0.55pc{
\\
r}
&
\\
\hline

\xymatrix@=1.5pc{
&*{\bullet} \ar@{-}[dl]\ar@{-}[d]\ar@{-}[r]& *{\bullet} \ar@{-}[d]\ar@{-}[dr] \\
*{\bullet} \ar@{-}[dr] &*{\bullet} \ar@{-}[d] & *{\bullet} \ar@{-}[d] &*{\bullet} \ar@{-}[dl] \\
&*{\bullet} \ar@{-}[r]& *{\bullet}  \\}

&
\xymatrix@=1.1pc{
\\
1}

&
\xymatrix@=0.98pc{
\\
{r}/{2}}
&
\xymatrix@=1.25pc{
\\
r}
&
\\
\hline

\end{tabular}

\vspace{1mm}
\

\end{center}

\begin{exa}[the special fibre is irreducible and has a single node]\label{exa:1irr}
This case is well know. The curve $C_K$ is smooth and its minimal regular model
is a stable curve whose special fibre contains a single node, which is nonseparating.
The group $\Pic_{C_K}[r]$ is ramified over
$K$ and, therefore, its N\'eron model over $R$ is not finite.

After base change to $K_d=K[\wt \pi]/(\wt \pi^d-\pi)$ we obtain a group
with finite N\'eron model if and only if $d$ is a multiple of $r$.
This is consistent with Theorem \ref{thm:group} and equation \eqref{eq:m1m2}, which yields
$m_1=1$, $m_2=r$, and $m_3=r$.
\end{exa}

\begin{exa}[the special fibre has two irreducible components with two nodes]\label{exa:2irr}
In this case the minimal regular model of $C_K$ is a stable curve whose special fibre
has two irreducible components and two (nonseparating) nodes.
On the one hand, this is a case where
the classical N\'eron model of the group of $2$-torsion line bundles
is finite as we illustrate hereafter.
On the other hand, the finite N\'eron model does not represent the group of
$2$-torsion line bundles of any
semistable reduction of $C_K$. Indeed this is a consequence of Remark \ref{rem:seqgl}:
as soon as $C_K$ admits a semistable reduction containing nonseparating nodes, the functor of
$r$-torsion bundles on any semistable reduction is not finite.

Let us focus on the finiteness of the N\'eron model or, equivalently, on the
proof of $m_2=1$ for the group $\Pic_{C_K}[2]$.
The dual graph is the second graph of the table above.
There is only one nontrivial
circuit in the graph and it has two edges; we have $c=2$.
By \eqref{eq:m1m2}, we have $m_2=1$; i.e. $N(\Pic_{C_K}[2])$ is finite.
In general, for $\Pic_{C_K}[r]$, we have $m_1=1$, $m_2=r/2$, and $m_3=r$.

\end{exa}
\begin{rem}
The previous example  shows that $m_2$ and $m_3$ are not equal in general.
This means that a finite  N\'eron $d$-model need not
represent the $r$-torsion line bundles on a twisted semistable reduction on $S[d]$.
\end{rem}
\begin{exa}[the dual graph of the special fibre is a circuit]
The above example can be generalised without change to the case where the
dual graph of the special fibre is a circuit having $h$ nodes and $h$ edges. In the
table we consider the case $h=4$.
If $r$ is a multiple of $h$, then we get $m_1=1$, $m_2=r/h$, and $m_3=r$.
For $h=r$ this provides a class of examples satisfying the necessary
and sufficient condition of Corollary \ref{cor:group}
for the finiteness of the N\'eron model of $\Pic_{C_K}[r]$.
\end{exa}

\begin{exa}\label{exa:3irr}
The fourth graph provides an example where,
for $r=4$, we have $m_1=1$, $m_2=2$, and $m_3=4$.
Since $m_2\ne 1$,  the N\'eron
(1-)model is not finite. Indeed, the condition of Corollary \ref{cor:group} is not satisfied
for $r=4$:
there are circuits sharing only two edges.
\end{exa}
\begin{rem}
For $r=4$, the two previous examples exhibit a situation where $C_K$  admits
a semistable minimal regular model
for which the dual graph of the special
fibre satisfies the following property.
The number of edges of the
circuits contained in the
dual graph of the special fibre is always a multiple of $r$.
Example \ref{exa:2irr} may lead us to
believe that this property implies the finiteness of the N\'eron model of
$\Pic_{C_K}[r]$.
Instead, Example \ref{exa:3irr} shows
that the N\'eron model of $\Pic_{C_K}[r]$ may not be finite.

The condition that the number of edges of
all circuits contained in the dual graph is always a multiple of $r$
is necessary, but not sufficient.
The right sufficient and necessary condition on the graph involves
the (signed) number  of edges shared by two circuits and is the
one stated in Corollary \ref{cor:group}.
\end{rem}
\begin{rem}\label{rem:lor}
We also point out that a
sufficient condition on the dual graph $\Gamma$ implying that
the N\'eron model of $\Pic_{C_K}[r]$ is finite is given by the following definition.
\end{rem}\begin{defn}\label{defn:rdivided}
A graph is \emph{$r$-divided} if
it is obtained by
taking another graph $\wt \Gamma$ and by dividing
all its edges in $r$ edges.
\end{defn}
\noindent Indeed, Proposition 2 of \cite{Lo} implies immediately
that
the N\'eron model of
$\Pic_{C_K}[r]$ is finite as soon as  $\Gamma$ is $r$-divided.
We illustrate with some examples that this condition is not necessary.

\begin{exa}[the dual graph of the special fibre need not be $r$-divided]\label{exa:wheels}
A simple example illustrating this fact is given by
taking
a graph with $2r+1$ edges, where the first $2r$ edges form
two disjoint circuits with $r$ edges and the last edge
joins two vertices lying on the first and the second circuit.
This graph is not $r$-divided (the fifth diagram
in the table provides a picture for $r=4$).

In fact, even if we restrict to  graphs which do not contain
separating edges, ``$r$-divided'' is stronger than
the condition stated in Corollary \ref{cor:group}. We show it in the
following example.
\end{exa}
\begin{exa}\label{exa:rdiv}
We consider the last graph of the table above.
All the circuits share an even number of edges, but the graph is not
$2$-divided.  Indeed, it contains two chains (the horizontal edges in the diagram)
 whose number of edges
is not a multiple of $2$.
\end{exa}

\section{Finite N\'eron  $d$-models of the torsor of $r$th roots}\label{sect:torsor}
We generalise the above results to the functor of $r$th roots $F_K^{1/r}$ of
a line bundle $F_K$ over the curve $C_K$ smooth over $K$. We assume that
the degree of $F_K$ is a multiple of $r$. The $K$-scheme $F_K^{1/r}$ is the
preimage of $F_K$ with respect to the $r$th tensor power in $\Pic_{C_K}$. It
has a natural structure of \'etale $K$-torsor under the kernel of the $r$th tensor
power: the $K$-group $\Pic_{C_K}[r]$.
\subsection{The finiteness of the N\'eron $d$-model of the torsor of $r$th roots.}
The following proposition exhibits a finite N\'eron $d$-models of $F_K^{1/r}$
that represents the functor of $r$th roots
over a twisted semistable reduction of $C_K$.
\begin{pro}\label{pro:torsor}
Under the  conditions of
Theorem \ref{thm:group},
we assume that $F_K$ is a line bundle on $C_K$,
whose degree is a multiple of $r$.
Recall that
$C_K$ has semistable reduction
on $S[d]$ if and only if $d$ is a multiple of $m_1$.

If $d$ is a multiple of $rm_1$, then the N\'eron $d$-model
of $F_K^{1/r}$ is finite
and there is a line bundle
$\sta F$ on a twisted semistable reduction $\stC$ of
$C_K$ on $S[d]$ extending $F_K\to C_K$
and satisfying
$$N_d(F_K^{1/r}) = \sta F^{1/r},$$
where $\sta F^{1/r}$ represents the functor
of $r$th roots of $\sta F$ on $\stC\to S[d]$.
In that case, the torsor structure of
$N_d(F_K^{1/r})$ under $N_d(\Pic_{C_K}[r])$
is the natural torsor structure of $\sta F^{1/r}$ under
the group of $r$-torsion line bundles on $\sta C$.
\end{pro}
\begin{proof}
Consider the semistable minimal regular model $C^{\rm reg}$ of $C_K$ on
$S[m_1]$. Since $C^{\rm reg}$ is regular, there exists a line bundle $F^{\rm reg}$ extending
$F_K$ on $C^{\rm reg}$.
For $d\in rm_1\ZZ$, over $S[d]$ there is a twisted curve $\stC$
whose nodes have stabiliser of order $r$ and whose coarse space fits in
\[
\xymatrix@R=0.1cm{
{\coa{\sta C}}\ar[rr]\ar[dd]&&C^{\rm reg}\ar[dd]\\
&\square&\\
S[d]\ar[rr]&&S[m_1]
}
\]
(we construct $\stC$ by iterating the construction of
Remark \ref{rem:modifsst} with $d=r$ at all nodes). Consider the pullback $\sta F$ of
$F^{\rm reg}$ on $\stC$ via the projection
to $C^{\rm reg}$ in the diagram above and via
 $\pi\colon \stC\to \coa{\stC}$.
Notice that $\sta F^{1/r}$ is a finite torsor on $S[d]$ by Theorem \ref{thm:cond}.
Then, for
$\wt R=R[\wt \pi]/(\wt \pi^d=\pi)$,
the pullback of $\sta F^{1/r}$ on $\Spec \wt R$ is the N\'eron model of $F_K^{1/r}$
by \cite[7.1/Thm.~1]{BLR}.
We conclude that $\sta F^{1/r}$
is the N\'eron $d$-model of $F_K^{1/r}$ by Proposition \ref{pro:relner}.
\end{proof}
\begin{rem}\label{rem:descentofF}
Note that the argument above shows that the line bundle $\sta F$
satisfying $N_{d}(F_K^{1/r})=\sta F^{1/r}$
can be chosen in such a way that it descends to a line bundle
over the semistable minimal regular model of $C_K$ on $S[m_1]$.
\end{rem}
\subsection{The finiteness criterion for $N(F_K^{1/r})$.}
We apply the above result to
determine when the
usual N\'eron model of $F_K^{1/r}$
is finite over $R$.
\begin{cor}\label{cor:torsor}
Let $r>2$ be prime to $\cha(k)$, and let $C_K$
be a smooth curve of genus $g\ge 2$. The N\'eron model of
$F_K^{1/r}$ is finite if and only if
$C_K$ admits a semistable minimal regular model
$C_R$ satisfying the following conditions:
\begin{enumerate}
\item The dual graph of the special fibre of $C_R$ verifies the condition
 of Corollary \ref{cor:group}.
\item There exists a line bundle $F_R$ on $C_R$
isomorphic to $F_K$ on the generic fibre and of degree $d\in r\ZZ$ on each irreducible component
of the special fibre.
\end{enumerate}\end{cor}
\begin{proof}
We note first
that the two conditions in the statement hold only if
$m_1$ equals $1$. Indeed, if the N\'eron model of the torsor is finite, then the N\'eron model of
$\Pic_{C_K}[r]$ is finite, hence $m_2=m_1=1$.

Using Proposition \ref{pro:torsor}, we regard the
N\'eron $r$-model of
$F_K^{1/r}$ as the torsor of $r$th roots
of a line bundle $\sta F$ extending
$F_K$ on a twisted reduction $\stC$ with stabilisers of
order $r$ at all nodes. Note that $\stC$ can be chosen in such a way that the
coarse space descends to the semistable minimal regular model $C^{\rm reg}_R$ of
$C_K$ on $S$. In Remark \ref{rem:descentofF}, we point
out that $\sta F$ descends to $C^{\rm reg}_R$.
Claiming that the N\'eron model of the torsor is finite
is equivalent to saying that $N_r(F_K^{1/r})$ also descends to $S$.
Denote by $\sta C_k$ the special fibre of $\stC$,
whose coarse space is the special fibre of
$C^{\rm reg}_R$; then, the descent of
$N_r(F_K^{1/r})$ to $S$ is equivalent to the claim that
the automorphism
$\sta g$ defined in
\eqref{eq:g}
acts trivially on the $r$th roots of $\sta F\rest{\stC_k}$.

Recall that the pullback $\sta j^*$ via
the embedding $\sta j$ of the singular locus into $\stC_k$, can be
regarded as a homomorphism from $\Pic \stC$ to $(\pmmu_r)^E=C_1(\Gamma,\pmmu_r)$.
By Proposition \ref{pro:aut_lb} and Remark \ref{rem:torsor},
$\sta g$ acts trivially on $\sta F\rest{\stC_k}^{1/r}$ if and only if
$\sta j^*(\sta F\rest{\stC_k}^{1/r})$ lies
in the image of $\delta_r\colon C_0(\Gamma,\pmmu_r)\to C_1(\Gamma,\pmmu_r)$.
This is equivalent to two inclusions:
first we require that $\im\delta_r$ includes
$\sta j^*( \Pic_{\stC_k}[r])$, second we require that there exists one root of
 $\sta F\rest{\stC_k}$
 whose image via $\sta j^*$ is contained in
 $\im \delta_r$.
The first inclusion $\im\delta_r\supseteq \sta j^*( \Pic_{\stC_k}[r])$
amounts to require that
$C_R$ satisfies the condition (1) of the statement above.
The second inclusion simply means that there exists an $r$th root $\sta L$
 of $\sta F$ on $\stC_k$ such that
$\im(\partial_r\circ \delta_r)$ contains
$\partial_r(\sta j^* \sta L)$.
A direct calculation shows that $\partial_r(\sta j^*\sta L)$
is the multiindex of degrees modulo $r$
of $\sta L^{\otimes r}\cong \sta F$
on each irreducible component of $\sta C_k$.
Therefore such a multiindex lies  in $\im(\partial_r\circ \delta_r)$
modulo $r$, which is, by Remark \ref{rem:FonCreg},
what needed to be shown.
\end{proof}

\begin{rem}For $r=2$, Remark \ref{rem:r=2} extends word for word. The criterion holds
in the case  $r=2$ if we assume that the minimal regular model $C_R$ of $C_K$ is semistable.
More explicitly, for any line bundle $F_K$ of even degree on $C_K$, we have $N(F_K^{1/2})$ if and only if
\begin{enumerate}
\item the circuits in the dual graph of the special fibre of $C_R$ always share an even number of edges.
\item there exists a line bundle $F_R$ on $C_R$ isomorphic to $F_K$ on the generic fibre of $C_R$ and
of even degree on each irreducible component of the special fibre of $C_R$.
\end{enumerate}
\end{rem}

\begin{exa}[square roots of a line bundle of multidegree $(1,-1)$.]\label{exa:banana}
We fix $r=2$, and we assume that the minimal regular model of $C_K$ is a
stable curve $C_R$,
whose special fibre $C_k$ has
two connected components $C_{1,k}$ and $C_{2,k}$ and two nodes.
%
%

Let us consider a line bundle $F_K$ of even degree over $C_K$.
Since $C_R$ is regular, it is the minimal regular model of $C_K$.
Furthermore, the study of the dual graph of the special fibre already
carried out in Example \ref{exa:2irr} implies that condition
(1) is  satisfied.
The finiteness of $N(F^{1/2}_K)$ only depends on the existence of a line bundle
$F_R$ extending $F_K$ with even degrees $d_1$ and $d_2$ on $C_{1,k}$ and $C_{2,k}$.

In fact, for any line bundle $F_R$ extending $F_K$ on $C_R$,
the degrees $d_1$ and $d_2$ have the same parity: $d_1+d_2\in 2\ZZ$.
If both degrees are even, then the condition (2) of
Corollary \ref{cor:torsor} is satisfied
and $N(F_K^{1/2})$ is finite. If both degrees are odd, then
the analysis carried out in
\S\eqref{sect:multideg} shows immediately that
there is no line bundle $F_R$ on $C_R$ extending $F_K$
with even degrees $d_1$ and $d_2$ on the special fibre.
Then, $N(F_K^{1/2})$ is not finite.

On the other hand, even in the case $d_1,d_2\in 1+2\ZZ$,
Proposition \ref{pro:torsor} says that the N\'eron $2$-model is finite.
Indeed its special fibre  represents $2^{2g}$
square roots on a twisted curve $\stC_k$ with stabiliser of order $2$ on each node.
The degree of each square root on the two irreducible components
of $\stC_k$ belongs to $1/2+\ZZ$.
\end{exa}

\section{The N\'eron model of $\Pic_{C_K}[r]$ without the theory of twisted curves}
\label{sect:without}
As stated in the introduction, this section is self contained
and can be read straight after $\S1$ assuming some
well known facts on semistable curves
recalled in \S2.1-5. We consider the question
of the finiteness of $\Pic_{C_K}[r]$ under the same conditions
of \S\eqref{sect:question}.

\subsection{The group of connected components of the special fibre of the N\'eron model.}
We consider the subfunctor $\Pic^0_{C_K}$ of
line bundles of degree 0.
It is represented by a proper, geometrically connected scheme of finite type
over $K$: the jacobian variety
of $C_K$. Let $N(\Pic^0_{C_K})$ be its N\'eron model over $R$.
In general, its special fibre
is not geometrically connected and
the connected components form a finite abelian group.

If $C_K$ admits a
semistable minimal regular model $C_R$, then
the group of components
of the special fibre of $N(\Pic_{C_K}^0)$ admits a simple description,
which is due to Raynaud \cite{Re_Pic}.
Indeed, we can realise this group by applying to
the dual graph $\Gamma$ of the special fibre of $C_R$ a
 construction
which associates to any connected graph a finite group.

\begin{defn}\label{defn:dcg}
Let $\Gamma$ be an oriented connected graph.
The group
$\Phi_\Gamma$ is an invariant of the
graph $\Gamma$ given by
\begin{equation}\label{eq:firstdefgroup}
\Phi_\Gamma= \frac{\partial({C_1}(\Gamma,\ZZ))}{\partial (\delta ({C_0}(\Gamma,\ZZ)))}.\end{equation}
\end{defn}

\begin{rem}\label{rem:dcgwithM}
Since $\Gamma$ is connected, the homomorphism $\partial \colon
C_1\to C_0$ and the augmentation homomorphism
$\varepsilon \colon C_0\to \ZZ$ form an exact sequence $C_1\to C_0\to \ZZ$.
We have
\begin{equation}\label{eq:kerim}
\Phi_\Gamma= \frac{\ker \varepsilon}{\im M},\end{equation}
where $M$ is the endomorphism of $C_0$ defined by $-\partial \circ \delta$.
(This endomorphism has a geometric interpretation, see \eqref{sect:multideg} and
in particular Remark \ref{rem:Mdelta}.)

In fact $\Phi_\Gamma$ is  a finite abelian group, \cite[9/10]{BLR}.
Furthermore,
write $h$ for the number of vertices of $V$,
the Smith normal form of $M$ is given by
${\rm Diag}(n_1, \dots , n_{h-1},0)$, with $n_i\ne 0$ all $i$.
In this way $\Phi_\Gamma$ is isomorphic to $\bigoplus_i\ZZ/n_i\ZZ$.
\end{rem}

\subsection{Definitions of $\Phi_\Gamma$ in the literature.}\label{sect:literature}
 This group appeared
 first in 1970, with the definition \eqref{eq:kerim},
in Raynaud's paper \cite[\S8.1.2, p. 64]{Re_Pic}. There, we consider
a curve $C_K$ whose minimal regular model $C_R$ is semistable.
The graph $\Gamma$ is the
dual graph of the special fibre
of $C_R$, and
$\Phi_\Gamma$ is identified with the
group of connected components of the special fibre of the N\'eron model of $\Pic^0_{C_K}$.

It is worth mentioning that,
in this perspective, the elements of
this group can be regarded as classes of
multidegrees of line bundles over a
semistable curve; therefore Caporaso
refers to $\Phi_\Gamma$ as the ``degree class group'', \cite{Cap}.

The same group was studied from different viewpoints and with different names:
Berman \cite{Be86}, Lorenzini \cite{Lo89},
Dahr \cite{Da90} (sandpile group),
Bacher--de la Harpe--Nagnibeda \cite{BDN97} (Picard group),
Biggs \cite{Bi99} (critical group).

\subsection{Reformulating the definition of $\Phi_\Gamma$ and describing its $r$-torsion.}
In view of the proof of our criterion
for the finiteness of $N(\Pic_{C_K}[r])$
we give a new characterisation of $\Phi_\Gamma$.

\begin{lem}
The following diagram is commutative and all
sequences are exact \begin{equation}\label{diag:exa}
\xymatrix@R=0.4cm{
         & 0                                           & 0               & 0                        & \\
0\ar[r]  & {H_1}/({H_1\cap \delta({C_0})})\ar[r]\ar[u] & H^1\ar[r]\ar[u] & \Phi_\Gamma\ar[r]\ar[u] & 0\\
0\ar[r]  & {H_1}  \ar[r]\ar[u] & {C_1}\ar[r]^{\partial}\ar[u]&\partial({C_1})\ar[r]\ar[u] & 0\\
0\ar[r]  & H_1\cap \delta({C_0})\ar[r]\ar[u] & \delta({C_0})\ar[r]^{\partial}\ar[u] & \partial(\delta ({C_0}))\ar[r]\ar[u] & 0\\
         & 0\ar[u] & 0\ar[u] & 0\ar[u] & ,\\
}
\end{equation}
where  $H_1$ and $H^1$ denote  $H_1(\Gamma,\ZZ)$ and $H^1(\Gamma,\ZZ)$,
and ${C_0}$ and ${C_1}$ denote  ${C_0}(\Gamma,\ZZ)$ and
${C_1}(\Gamma, \ZZ)$.
In this way, we have
\begin{equation}\label{eq:H1}\Phi_\Gamma\cong \frac{H^1}{{H_1}/({H_1\cap \delta({C_0})})}.
\end{equation}\qed
\end{lem}

Since $\Phi_\Gamma$ is a finite abelian group
its $r$-torsion $\Phi_\Gamma[r]$ is given by
$\Phi_\Gamma\otimes \ZZ/r\ZZ$. By tensoring both sides of
\eqref{eq:firstdefgroup}, we get
$$\Phi_\Gamma[r]\cong \frac{\partial_r({C_1}(\Gamma,\ZZ/r\ZZ))}
{\partial_r(\delta_r({C_0}(\Gamma,\ZZ/r\ZZ)))}.$$
Write the analogue of diagram \eqref{diag:exa} with ${C_0}(\Gamma,\ZZ/r\ZZ)$,
${C_1}(\Gamma,\ZZ/r\ZZ)$, $\partial_r$, and $\delta_r$ instead of
${C_0}(\Gamma,\ZZ)$, ${C_1}(\Gamma,\ZZ)$, $\partial$, and $\delta$.
The diagram is still exact, provided that we replace $H_1$ and $H^1$ by
$H_1(\Gamma,\ZZ/r\ZZ)=\ker (\delta_r)$ and $H^1(\Gamma,\ZZ/r\ZZ)=\coker \delta_r$.
In this way, we get the following exact sequence.
\begin{lem}\label{lem:rgroup} We have
\begin{multline}
0\rightarrow {{H_1(\Gamma,\ZZ/r\ZZ)}/({H_1(\Gamma,\ZZ/r\ZZ)\cap \delta_r({C_0(\Gamma,\ZZ/r\ZZ)})})}
\longrightarrow H^1(\Gamma,\ZZ/r\ZZ) \longrightarrow  \Phi_\Gamma[r]\rightarrow  0
\end{multline}\qed
\end{lem}
Recall that $H^1(\Gamma,\ZZ/r\ZZ)$ is isomorphic to $(\ZZ/r\ZZ)^{\oplus b_1(\Gamma)}$.
Therefore, $\Phi_\Gamma[r]$ is a quotient of
a group of order $r^{2g}$.

\subsection{The group $\Phi_\Gamma$ and the finiteness of the
N\'eron model of $\Pic_{C_K}[r]$.}
Applying Raynaud's description of the N\'eron model of
the jacobian variety, the condition that $N(\Pic_{C_K}[r])$
is finite can be characterised in terms of the group $\Phi_\Gamma$.
\begin{pro}\label{pro:viane}
Let $r>2$ be prime to $\cha(k)$. The N\'eron model of $\Pic_{C_K}[r]$ is finite if and only if
the minimal regular model of $C_K$ over $R$ is semistable
and the dual graph $\Gamma$ of its  special fibre satisfies
\begin{equation}\label{eq:gbetti}
\Phi_\Gamma[r]\cong(\ZZ/r\ZZ)^{\oplus b_1(\Gamma)}.\end{equation}
\end{pro}
\begin{proof}
Since $r>2$, the fact that the
N\'eron model $N(\Pic_{C_K}[r])$ is finite implies that
the minimal regular model $C_R$ of $C_K$ over $R$
is semistable (see \cite[5.17]{De}).
Therefore, we can take for granted the existence of
a semistable regular reduction $C_R$.

Then the $r$-torsion of the special fibre of
$N(\Pic^0_{C_K})$ is the special fibre
of the $r$-torsion of $N(\Pic^0_{C_K})$, which is the
N\'eron model of $\Pic^0_{C_K}[r]=\Pic_{C_K}[r]$.
In this way, the
N\'eron model $N(\Pic_{C_K}[r])$ is finite, if and only if
the special fibre  of the N\'eron model $N(\Pic_{C_K}[r])$
consists of $r^{2g}$ points.
By \cite{Re_Pic} and \cite[9]{BLR}, the special fibre of
$N(\Pic_{C_K}[r])$ is  an extension of
$\Phi_\Gamma[r]$ by the $r$-torsion of the Picard group
of the special fibre $C_k$ of $C_R$.
Now,
$\Pic_{C_k}[r]$ consists of $r^{2g-b_1(\Gamma)}$ points.
Therefore, the number of points of the special fibre of
$N(\Pic_{C_K}[r])$ is $r^{2g}$ if and only if
the group $\Phi_\Gamma[r]$ has order $r^{b_1}$. This fact,
by Lemma \eqref{lem:rgroup},
 is equivalent  to $\Phi_\Gamma[r]\cong (\ZZ/r\ZZ)^{\oplus b_1(\Gamma)}$.
\end{proof}

\subsection{When is the $r$-torsion subgroup of $\Phi_\Gamma$ isomorphic to $(\ZZ/r\ZZ)^{\oplus b_1}$?}
We obtain again the answer of Corollary \ref{cor:group}.
\begin{pro}\label{pro:tor} For any positive $r$, the $r$-torsion subgroup of
the group $\Phi_\Gamma$ is isomorphic to $(\ZZ/r\ZZ)^{\oplus b_1(\Gamma)}$ if and only if
the (signed) number of edges common to any two  circuits is always a multiple of $r$.
\end{pro}
\begin{proof}
Recall the exact sequence of Lemma \ref{lem:rgroup} and
the fact that $H^1(\Gamma,\ZZ/r\ZZ)$ is isomorphic to $(\ZZ/r\ZZ)^{\oplus b_1(\Gamma)}$.
Then, $\Phi_\Gamma[r]$ is isomorphic to $(\ZZ/r\ZZ)^{\oplus b_1(\Gamma)}$
if and only if
${H_1(\Gamma,\ZZ/r\ZZ)}={H_1(\Gamma,\ZZ/r\ZZ)\cap \delta({C_0(\Gamma,\ZZ/r\ZZ)})}.$
This yields $\ker (\partial_r)\subseteq \im(\delta_r)$.
Lemma \ref{lem:circuits} proven in Section \ref{sect:circuits}
implies the statement of the theorem.
\end{proof}

\section{The proof of the combinatorial lemma on circuits} \label{sect:circuits}
We give a proof of a combinatorial lemma used
in Proposition \ref{pro:m2m1},
and in Proposition \ref{pro:tor}.

\begin{lem}\label{lem:circuits}
For any positive integer
$r$ and for any connected graph $\Gamma$ the following conditions are equivalent:
\begin{enumerate}
\item $\ker(\partial_r) \subseteq \im (\delta_r)$;
\item the  number  of edges common to any two  circuits is always a multiple of $r$.
\end{enumerate}\end{lem}
\begin{proof}[Proof of Lemma \ref{lem:circuits}.]
Using $(2)$, we show that  $\ker(\partial_r)$
is included in
$\im(\delta_r)$.
It is enough to prove this claim for any circuit $C=\sum _{i=0}^{n-1} [e_i]$
passing through $n$ vertices $v_0, v_1, \dots, v_{n-1}$:
for $0\le i< n-1$,
$v_i$ to $v_{i+1}$ are the tail and the tip of $e_i$,
whereas $v_{n-1}$ and $v_{0}$ are the tail and the tip of
$e_{n-1}$. We have
$C\in\ker(\partial_r)$, and we want to show that there exists
$A=\sum_{v\in V} a(v)[v]$
such that $\delta_r(A)=C$.

For an arbitrary vertex $v$, we define $a(v)$ as follows: if there exists a
simple path joining
$v$ to $v_i$ with all vertices outside $C$ apart from $v_i$, then we set $a(v)\equiv i \mod r$.
Note that it may happen that a vertex  $v$ off $C$ can be connected
in the way described above to both
$v_i$ and $v_j$ with $i<j$ via
two simple paths $S_i$ and $S_j$. Then,
we can choose the last vertex $l$
common to $S_i$ and $S_j$, and
define a circuit $\wt C$ joining $l$ to $v_i$ via $S_i$,
$v_i$ to $v_j$ via $C$, and, finally, $v_j$ to $l$ via
$-S_j$. The condition $(2)$ applied to $\wt C$ and $C$ yields
$j-i\equiv 0 \mod r$.
It is easy to check that, for $A=\sum_{v\in V} a(v)[v]$, we have $\delta_r(A)=C$ by construction.

Conversely, we consider two  circuits $C$ and $\wt C$
and we count the common edges (with sign).
We write $\wt C$ as $\sum_{j=1}^{m-1}[f_j]$,
the circuit whose vertices are $w_0$, \dots, $w_{m-1}$, whose $j$th edge
$f_j$ is oriented from $w_j$ to $w_{j+1}$ for $0\le j< m-1$, and whose
edge $f_{m-1}$ is oriented from $w_{m-1}$ to
$w_{0}$. By $\ker (\partial_r)\subseteq \im(\delta_r)$
there exists $A=\sum _{v\in V} a(v)[v]\in C_0(\Gamma,\ZZ/r\ZZ)$
such that $\delta_r(A)=\wt C$. Using $A$, for any $j\in \{0, \dots, m-1\}$,
we associate to
each vertex $w_j$ an index $a(w_j)\in \ZZ/r\ZZ$,
which we denote by $a_j\in \ZZ/r\ZZ$.
The condition $\partial_r(A)=\wt C$ implies
\begin{align}\label{eq:x-y}
m\in r\ZZ  \quad \quad \text{and} \quad \quad a_x-a_y=x-y  \quad \quad \text{in }\ZZ/r\ZZ,
\end{align}
if $x$ and $y$ are values in $\{0, \dots, m-1\}$.
%
We also observe that, by definition,
\begin{equation}\label{eq:getrid}
a(v_+)=a(v_-) \quad \quad \text{in }\ZZ/r\ZZ\end{equation}
if $v_+$ and $v_-$ are
the tip and the tail  of an  edge not included in $\wt C$.

We assume that $C$ and $\wt C$ have some edges in common, but not all
(otherwise the claim follows from $m\in r\ZZ$ in \eqref{eq:x-y}).
We write $C=\sum _{i=0}^n [e_i]$ for the circuit $C$
passing through the vertices $v_0, \dots , v_{n-1}$. We choose the notation so that,
for $0\le i< n-1$,
the edge $e_i$ is oriented from $v_i$ to $v_{i+1}$, whereas
the edge
$e_{n-1}$ is oriented from $v_{n-1}$ to $v_{0}$
and \emph{does not lie on $\wt C$}. By \eqref{eq:getrid},
this implies
\begin{equation}\label{eq:getrid0}
a(v_{n-1})=a(v_0) \quad \quad \text{in }\ZZ/r\ZZ.
\end{equation}

Consider the indices $0,1$, \dots, $n-1$ of the vertices of $C$.
Define a non decreasing
sequence of integers $0\le t_1\le  \dots\le  t_{2h}<n$ with $h\in \NN$
so that $i$ and $i+1$ are contained in
one of the  intervals
$[t_1, t_2], [t_3,t_4],
\dots, [t_{2h-1}, t_{2h}]$
if and only if there exists an edge of
$\wt C$ joining $v_i$ and $v_{i+1}$ (in any direction).
Note that we can insert in the sequence two
values $t_{2j-1}$ and $t_{2j}$ with $t_{2j-1}=t_{2j}$
for any $j\in \NN$,
as long as we reparametrise and we respect the monotony.
In this way,
we can adjust the choice of the values $\{t_i\}_{0<i\le 2h}$ so that $h$ is even
and for each even (resp. odd) $k$
satisfying $0< k \le h$,  the simple path
$\sum _{t_{2k-1} \le i < t_{2k}} e_{i}$
joining $v_{t_{2k-1}}$ to $v_{t_{2k}}$
lies on $\wt C$ and $C$ with the same (resp. the opposite) orientation.
Now, \eqref{eq:x-y} implies
\begin{align}\label{eq:trick}
a(v_{t_{2k}})-a(v_{t_{2k-1}})=(-1)^k (t_{2k}-t_{2k-1}) \quad \quad \text{in }\ZZ/r\ZZ;
\end{align}
Clearly, the following relation hold in $\ZZ/r\ZZ$
$$
0= (a(v_0) -a (v_{n-1}))+\textsum_{0\le j<n-1} (a(v_{j+1})-a(v_{j})).$$
By \eqref{eq:getrid0}, the first summand on the right hand side vanishes.
In fact, using \eqref{eq:getrid}, we get
$$0=\textsum_{0<k\le h} (a(v_{t_{2k}})-a(v_{t_{2k-1}})).$$
By \eqref{eq:trick}, we have
\begin{align*}
0=\textsum_{0<k\le h} (-1)^k(t_{2k}-t_{2k-1}).
\end{align*}
Note that the last equation in $\ZZ/r\ZZ$ is the claim $(2)$.
\end{proof}

\small
{
\font\cc=cmcsc10
\def\xx#1{\par\noindent\hangindent=1cm\hangafter=1{\cc {#1},}}

\vskip.5cm
\noindent\textsc{A.~Chiodo,}
Laboratoire Painlev\'e (Lille 1), UMR du CNRS 8524,
       59655 Villeneuve d'Ascq cedex, France. \\
       \noindent email: \texttt{chiodo@unice.fr}}
\end{document}